\documentclass[aap,MSNbibl,seceqn,citesort,dvips]{arximspdf}

\doi{10.1214/11-AAP818} 
\volume{22}
\issue{5}
\pubyear{2012}
\firstpage{1928}
\lastpage{1961}

\makeatletter

\newtheorem{theorem}{Theorem}[section]
\newtheorem{lemma}[theorem]{Lemma}
\newtheorem{proposition}[theorem]{Proposition}

\newproclaim{remark}{Remark}

\newcommand{\GA}{D}
\newcommand{\YY}{Z}
\newcommand{\QQ}{{\mathcal{S}}}
\newcommand{\be}{\varpi}

\newcommand{\tphi}{w}
\newcommand{\ZZ}{\mathcal{Z}}
\newcommand{\W}{{l}}

\newcommand{\CF}{\mathcal{F}}
\newcommand{\CB}{\mathcal{B}}
\newcommand{\CL}{\mathcal{L}}
\newcommand{\CI}{\mathcal{J}}
\newcommand{\CG}{\mathcal{G}}
\newcommand{\CS}{\mathcal{S}}

\renewcommand{\b}{\beta}
\newcommand{\eps}{\varepsilon}
\newcommand{\Th}{\Theta}
\newcommand{\s}{\sigma}

\newcommand{\J}{\mathbb{I}}
\newcommand{\N}{\mathbb{N}}
\newcommand{\E}{\mathbb{E}}
\newcommand{\RE}{\mathbb{R}}
\newcommand{\setN}{\mathbb{N}}
\newcommand{\T}{\mathbb{T}}
\newcommand{\U}{\mathbb{U}}

\newcommand{\DXi}{\frac{\partial}{\partial\xi}}

\newcommand{\DT}{\frac{\partial}{\partial t}}

\makeatother

\begin{document}
\begin{frontmatter}

\title{Self-similar solutions in one-dimensional kinetic models: A probabilistic view}
\runtitle{Self-similar solutions in one-dimensional kinetic models}

\begin{aug}
\author[A]{\fnms{Federico} \snm{Bassetti}\ead[label=e1]{federico.bassetti@unipv.it}}
\and
\author[B]{\fnms{Lucia} \snm{Ladelli}\corref{}\ead[label=e2]{lucia.ladelli@polimi.it}}
\runauthor{F. Bassetti and L. Ladelli}
\affiliation{Universit\`a degli Studi di Pavia
and Politecnico di Milano}
\address[A]{Dipartimento di Matematica \\
Universit\`a degli Studi di Pavia \\
via Ferrata 1\\
27100, Pavia\\
Italy\\
\printead{e1}}
\address[B]{Dipartimento di Matematica \\
Politecnico di Milano \\
piazza Leonardo da Vinci 32\\
20133, Milano\\
Italy \\
\printead{e2}} 
\end{aug}

\received{\smonth{4} \syear{2010}}
\revised{\smonth{9} \syear{2011}}

%
\begin{abstract}
This paper deals with a class of Boltzmann equations on the real line,
extensions of the well-known Kac caricature. A distinguishing feature
of the corresponding equations is that therein, the collision gain
operators are defined by $N$-linear smoothing transformations. These
kind of problems have been studied, from an essentially analytic
viewpoint, in a recent paper by Bobylev, Cercignani and Gamba
[\textit{Comm. Math. Phys.} \textbf{291} (2009) 599--644]. Instead, the
present work rests exclusively on probabilistic methods, based on
techniques pertaining to the classical central limit problem and to the
so-called fixed-point equations for probability distributions. An
advantage of resorting to methods from the probability theory is that
the same results---relative to self-similar solutions---as those
obtained by Bobylev, Cercignani and Gamba, are here deduced under
weaker conditions. In particular, it is shown how convergence to a
self-similar solution depends on the belonging of the initial datum to
the domain of attraction of a specific stable distribution. Moreover,
some results on the speed of convergence are given in terms of
Kantorovich--Wasserstein and Zolotarev distances between probability
measures.
\end{abstract}

%
\begin{keyword}[class=AMS]
\kwd[Primary ]{60F05}
\kwd[; secondary ]{82C40}.
\end{keyword}
\begin{keyword}
\kwd{Central limit theorem}
\kwd{domain of normal attraction}
\kwd{stable law}
\kwd{Kac model}
\kwd{smoothing transformations}
\kwd{marked recursive $N$-ary random trees}
\kwd{self-similar solution}.
\end{keyword}

\end{frontmatter}

\section{Introduction}\label{intro}

In this paper we consider a kinetic-type evolution equation, introduced
and studied in \cite{CeGaBo}, which includes some well-known
one-dimensional Maxwell models. If $\phi(t,\xi)$ denotes the
Fourier--Stieltjes transform
\[
\phi(t,\xi):=\int_\RE e^{i\xi v}\rho_t(dv) \qquad(\xi\in\RE)
\]
of a time-dependent probability measure $\rho_t$ on the real line $\RE$,
the equation under interest is
%
\begin{equation}
\label{eqboltzivp}\quad
\cases{
\displaystyle \DT\phi(t,\xi)+\phi(t,\xi) =
\hat{Q}(\phi(t,\cdot),\ldots, \phi(t,\cdot)) (\xi)
\qquad(t>0, \xi\in\RE),\vspace*{2pt}\cr
\phi(0,\xi)=\phi_0(\xi),}
\end{equation}
where, given $N$ characteristic functions $\phi_1,\ldots,\phi_N$,
%
\begin{equation}
\label{eq2}
\hat{Q}(\phi_1,\ldots,\phi_N)(\xi)
:= \E[\phi_1(A_1\xi) \cdots\phi_N(A_N\xi)] \qquad(\xi
\in\RE).
\end{equation}
The expectation $\E$ in (\ref{eq2}) is taken with respect to the
distribution of a given vector $A=(A_1,\ldots,A_N)$ of positive
real-valued random variables
defined on a probability space $(\Omega,\CF,P)$.
The initial condition $\phi_0$ is a characteristic function of a prescribed
real random variable $X_0$ with distribution function $F_0(x)$.

Notice that different equations for probability
dynamics considered in literature are special cases of (\ref{eqboltzivp}):
the one-dimensional Kac caricature \cite{Kac}, some one-dimensional
dissipative Maxwell models
\cite{Ben-Avraham,PareschiToscani,PulvirentiToscani}, some mean
conservative models used to
describe economical dynamics (see, e.g., \cite{MatthesToscani,Pat}), some
models for mixture of Maxwell gases \cite{BoGa2006}.
In addition, using results in \cite{CeGaBo,CeGaBoBis}, it can be shown
that the isotropic solutions of the multidimensional inelastic
Boltzmann equation \cite{BoCe}
are functions of one-dimensional Fourier--Stieltjes transforms which
are solutions
of (\ref{eqboltzivp}) for a suitable choice of $(A_1,\ldots,A_N)$
and $\phi_0$.
Finally, we recall that this kind of kinetic equations
describes the evolution of the limit (in a suitable sense)
of a pure jump Markov process, representing $K$ interacting particles,
when $K$ diverges; see, for instance, \cite{MeGr}.

For simplicity of notation in the rest of the paper we write $\hat
Q(\phi)$ instead of $\hat Q(\phi,\ldots,\phi)$.

The aim of this paper is to study the
asymptotic behavior of
the solution~$\phi$ of (\ref{eqboltzivp}) as $t \to+\infty$.

One can distinguish two different situations:
\begin{itemize}
\item the solution $\phi(t,\xi)$ converges, as $t \to+\infty$, to a
stationary solution, that is, a~characteristic function
$\phi_\infty$ such that
%
\begin{equation}\label{stazeq}
\phi_\infty= \hat Q(\phi_\infty);
\end{equation}
\item there exists $\mu^*$ (depending on the initial
condition $\phi_0$) such that the rescaled solution
%
\begin{equation}\label{rescaledsol}
\tphi(t,\xi):=\phi(t,e^{-\mu^*t}\xi)
\end{equation}
converges as $t
\to+\infty$ to a nondegenerate limit.
\end{itemize}
To understand the nature of
this limit, let us observe that the re-scaled solution, $\tphi$,
satisfies the following new equation:
%
\begin{equation}
\label{eq2bis}
\cases{
\displaystyle \DT\tphi(t,\xi)+ \mu^* \xi\,\DXi\tphi(t,\xi)+\tphi(t,\xi) =
\hat{Q} (\tphi(t,\cdot))(\xi),\vspace*{2pt}\cr
\tphi(0,\xi)=\phi_0(\xi).}
\end{equation}
When $\mu^*=0$ equation (\ref{eq2bis}) reduces to (\ref{eqboltzivp})
and, clearly,
$\tphi$ is simply $\phi$.
The stationary equation associated to (\ref{eq2bis})
is, for every $\mu^*$,
\[
\mu^* \xi\,\DXi\tphi_\infty(\xi)+\tphi_\infty(\xi) =
\hat{Q}(\tphi_\infty)(\xi),
\]
which can be re-written, after easy computations, as an integral
equation for a Fourier--Stieltjes transform
%
\begin{equation}\label{stazselfintegral}
\tphi_\infty(\xi)= \int_0^1 \hat{Q} (\tphi_\infty)(\tau
^{\mu^*} \xi
)\,d\tau.
\end{equation}

It is important to note that, if
a characteristic function $\tphi_\infty$ satisfies (\ref{stazselfintegral}),
then
\[
\phi(t,\xi):=\tphi_\infty( \exp\{\mu^*t \}\xi)
\]
satisfies the original Kac-like equation (\ref{eqboltzivp})
with initial condition $\phi_0(\xi)=\tphi_\infty(\xi)$.
Following \cite{CeGaBo}, we shall use the name \textit{self-similar
solution} for a solution $\tphi_\infty$
of (\ref{stazselfintegral}) (when it exists), although
the name self-similar solution is usually devoted to $\tphi_\infty
(\exp
\{\mu^*t \}\xi)$.

In terms of random variables, (\ref{stazselfintegral}) becomes
%
\begin{equation}
\label{eq3}
X\stackrel{\CL}=\Th^{\mu^*} \sum_{i=1}^N A_i X_i,
\end{equation}
where $(X,X_1,\ldots,X_N)$ are stochastically independent random variables
with the same characteristic function $\tphi_\infty$, $\Th$ is a random
variable with uniform distribution on $(0,1)$ and
$(X,X_1,\ldots,X_N)$, $\Th$ and $(A_1,\ldots,A_N)$
are stochastically independent. Moreover, $Z_1\stackrel{\CL}=Z_2$
means that
the random variables $Z_1$ and $Z_2$ have the same law.
$\hat Q$ is usually called smoothing transformation,
and equations of kind~(\ref{stazselfintegral}), or equivalently (\ref
{eq3}), are referred to as fixed point equations
for distributions.

In \cite{CeGaBo} Maxwell models of type (\ref{eqboltzivp}) are
considered from a very general point of view
and some key
properties that lead to the self-similar asymptotics are established
mainly by analytic techniques.
The goal of our paper is to study convergence to self-similar solutions
by means of probabilistic methods. Via a suitable probabilistic
representation of the solution of
(\ref{eqboltzivp})
we resort to \textit{central limit theorems} and
\textit{fixed point equations for distributions}. In this way we are able
to extend some
results presented in \cite{CeGaBo}. The main result we obtained is the
proof of long-time convergence of the rescaled solution
to a self-similar solution essentially under
the natural hypothesis that the initial condition belongs to the domain
of normal attraction of a stable distribution.
Our approach is a generalization of the methods developed in
\cite{BLM2008}, where only the convergence to  stationary solutions for
(\ref{eqboltzivp}) and (\ref{eq2}) with $N = 2$ has been studied.
We mention that in \cite{BaLaTo} a probabilistic approach has been used
to study the solutions
of a kinetic equation in which the collision gain operator is a
(bilinear) inhomogeneous smoothing transformation.

The general idea of representing solutions to Kac-like equations in a
probabilistic way
dates back at least to \cite{McKean1966};
this approach has been fully formalized and employed in the derivation
of various results in
the last decade; see, for example, \cite
{CarlenCarvalhoGabetta2000,GabettaRegazziniCLT}.
For the original Kac equation, probabilistic methods have been used in
many papers;
see \cite{RegazUmi} for a review.

The paper is organized as follows. Section~\ref{SecMain} contains the
statements of our main theorems.
In Section~\ref{Srandtreeandprobrep} we derive the stochastic
representation of solutions to
(\ref{eqboltzivp}). Section~\ref{Smartingaleetal}
contains the statements of some intermediate results concerning sums of
random variables
indexed by random $N$-ary recursive trees.
All proofs are completed in Section~\ref{Sproofs}.

\section{Main results}\label{SecMain}

From now on we assume that
$A_i$ are nonnegative random variables
such that
%
\begin{equation}\label{H0-1}
P \Biggl\{ \sum_{i=1}^N \J\{A_i>0\} \in\{0,1\} \Biggr\}<1,\qquad  \E
\Biggl[\sum_{i=1}^N \J\{A_i >0 \} \Biggr]>1.
\end{equation}

In the theorems below the initial condition $F_0$ will satisfy one of
the following
hypotheses~\ref{hypoHgama}, where $\gamma$ belongs to $(0,2]$:

{\renewcommand\thelonglist{(${\mathbf H_{1}}$)}
\renewcommand\labellonglist{\thelonglist}
\begin{longlist}
\item\label{hypoH1}
\textit{either} (a) $\int_\RE|v| \,dF_0(v) <+\infty$ \textit{and}
$m_0=\int_\RE v\,dF_0(v)$
\textit{or} (b) $F_0$ \textit{is a symmetric distribution function and
satisfies the condition}
%
\begin{equation}
\label{stabledomain2}\quad
\lim_{x \to+\infty} x \bigl(1-F_0(x)\bigr) =c_0^+<+\infty,\qquad
\lim_{x \to-\infty} |x| F_0(x) =c_0^+<+\infty
\end{equation}
\textit{with} $c_0^+>0$;
\end{longlist}}

{\renewcommand\thelonglist{(${\mathbf H_{2}}$)}
\renewcommand\labellonglist{\thelonglist}
\begin{longlist}
\item\label{hypoH2} 
$0<\s_0^2:=\int_\RE|v|^2 \,dF_0(v) <+\infty$ \textit{and}
$\int_\RE v\,dF_0(v)=0$;\\[2pt]
\textit{if} $\gamma\in(0,1) \cup(1,2)$
\end{longlist}}

{\renewcommand\thelonglist{(${\mathbf H_{\gamma}}$)}
\renewcommand\labellonglist{\thelonglist}
\begin{longlist}
\item\label{hypoHgama}
\textit{$F_0$ satisfies the condition}
%
\begin{equation}
\label{stabledomain}\qquad
\lim_{x \to+\infty} x^\gamma\bigl(1-F_0(x)\bigr) =c_0^+<+\infty,\qquad
\lim_{x \to-\infty} |x|^\gamma F_0(x) =c_0^-<+\infty
\end{equation}
\textit{with} $c_0^++ c_0^->0$
\textit{and}, \textit{in addition}, $\int_\RE v\,dF_0(v)=0$ \textit{if $\gamma>1$.}
\end{longlist}}

Accordingly, define
%
\begin{equation}
\label{chaSta}
\hat g_\gamma(\xi):=
\cases{
\displaystyle e^{ i m_0 \xi}, \hspace*{19.2pt}\qquad \mbox{if $\gamma=1$ and (a) of
\ref{hypoH1} holds,}\vspace*{2pt}\cr
\displaystyle e^{ - \pi c_0^+ |\xi| },\hspace*{5.9pt}\qquad \mbox{if $\gamma=1$ and
\textup{(b)} of~\ref{hypoH1} holds,}\vspace*{2pt}\cr
\displaystyle e^{ - \sigma_0^2 |\xi|^2/2 },\qquad \mbox{if $\gamma=2$ and
\ref{hypoH2} holds,}\vspace*{2pt}\cr
\displaystyle e^{ -k_0 |\xi|^\gamma
(1-i \eta_0 \tan({\pi\gamma}/{2} )\operatorname{sign}\xi) },
\cr
\hspace*{47.8pt}\qquad \mbox{if $\gamma\in(0,1) \cup(1,2)$
and~\ref{hypoHgama} holds,}}
\end{equation}
where
%
\begin{equation}
\label{constant}
k_0 = (c_0^{+}+c_0^{-}) \frac{\pi}{2\Gamma(\gamma)\sin(\pi\gamma/2)},
\qquad \eta_0 = \frac{c_0^{+}-c_0^{-} }{c_0^{+}+c_0^{-}}.
\end{equation}

Notice that, except when $\gamma=1$ and (a) of~\ref{hypoH1} holds,
$\hat g_\gamma$ is the Fourier--Stieltjes transform
of \textit{a centered stable law} of exponent $\gamma$ and
(\ref{stabledomain}) of~\ref{hypoHgama} is equivalent to saying that $F_0$
belongs to the domain of normal attraction of a $\gamma$-stable law
$g_\gamma$
with Fourier--Stieltjes transform $\hat g_\gamma$.
See, for example, Chapter 17 of \cite{fristedgray}.

It is worthwhile to recall that a distribution function $F_0$
belongs to the {domain of normal attraction} of
a stable law of exponent $\gamma$
if for any sequence of independent and identically distributed real-valued
random variables $(X_n)_{n \geq1}$
with common distribution function $F_0$,
there exists a sequence of real numbers $(c_n)_{n \geq1}$
such that the law of
$
{n^{-1/\gamma}} \sum_{i=1}^n X_i -c_n
$
converges weakly to a stable law of exponent $\gamma\in(0,2]$.

\subsection{Convergence to self-similar solutions}\label{SecMain1}
Our main result states that, under suitable assumptions, the rescaled
solution $w$, defined in (\ref{rescaledsol}), converges to
a mixture of centered stable characteristic functions.
{ The probability distribution of the mixing measure is
characterized as a particular (positive) solution of the fixed point equation
%
\begin{equation}\label{eqvRandomvariables}
\YY\stackrel{\CL}=\Th^{\QQ(\gamma)} \sum_{i=1}^N A_i^{\gamma}
\YY_i,
\end{equation}
where, as usual, $(\YY,\YY_1,\ldots,\YY_N)$ are i.i.d., $\Th$ is a
random variable with uniform distribution on $(0,1)$,
$(\YY,\YY_1,\ldots,\YY_N)$, $\Th$ and $(A_1,\ldots,A_N)$
are stochastically independent and
$\QQ\dvtx[0,\infty)\to[-1,\infty]$ is the convex function defined by
%
\begin{equation}
\label{eqthes}
\QQ(s):=\E\Biggl[\sum_{j=1}^N A_j^s \Biggr]-1
\end{equation}
with the convention that $0^0=0$.
In other words, }
the Fourier--Stieltjes transform of the mixing measure
will be characterized as a particular solution of
the integral equation
%
\begin{equation}\label{eqv}
v(\xi)=\int_0^1 \E\Biggl[\prod_{i=1}^N v\bigl(\xi A_i^\gamma\tau^{\QQ
(\gamma
)}\bigr) \Biggr] \,d\tau.
\end{equation}
Note that, thanks to (\ref{H0-1}), one has $0<\QQ(0)\leq N-1$; hence
if $\QQ(s)<+\infty$ for some $s$, then $\QQ(q)<+\infty$
for every $q$ in $(0,s)$.
We point out that there is a simple connection between
the function $s \mapsto\QQ(s)$, widely used in the probabilistic fixed
point literature, and
the so-called spectral function $s \mapsto{\mu}(s)$ introduced in
\cite
{CeGaBo}, more exactly
\[
{\mu}(s):=\frac{\QQ(s)}{s} \qquad(s >0).
\]
Finally, we observe that the function $\QQ(s)$ is such that
\[
\E\Biggl[\Th^{\QQ(s)} \sum_{i=1}^N A_i^{s}\Biggr]=1
\]
for every $s$, provided $\QQ(s)<+\infty$.
In the next proposition we collect some useful results concerning equation
(\ref{eqvRandomvariables}), or equivalently (\ref{eqv}).
%
\begin{proposition}\label{Pfixpoint} Fix $\gamma$ in $(0,2]$. Assume that
${\mu}(\delta)<{\mu}(\gamma)<+\infty$ for some $\delta>\gamma$. Then:
\begin{longlist}[(iii)]
\item[(i)] there is a unique probability distribution $\zeta_{\infty
,\gamma}$
on $\RE^+$, with
\[
\int_{\RE^+} z \zeta_{\infty,\gamma}(dz)=1,
\]
such that if
$\YY$ has law $\zeta_{\infty,\gamma}$, then it satisfies (\ref
{eqvRandomvariables}), or equivalently
$v_{\infty, \gamma}(\xi):=\int_{\RE^+} e^{i\xi z} \zeta_{\infty
,\gamma}(dz)$
is a solution of (\ref{eqv});
\item[(ii)] the equation
${\mu}(q)-{\mu}(\gamma)=0$ has at most one solution $q^*_\gamma\not
=\gamma
$, and we set, by convention, $q^*_\gamma:=+\infty$
if the unique solution is $q=\gamma$;
\item[(iii)] $\zeta_{\infty,\gamma}$ is degenerate if and only
$\sum_{i=1}^N A_i^\gamma=1$ almost surely. Moreover, if $P\{\sum
_{i=1}^N A_i^\gamma=1 \}<1$ and $p>\gamma$,
$\int_{\RE^+} z^{{p}/{\gamma}} \zeta_{\infty,\gamma
}(dz)<+\infty$
if and only if $p < q^*_\gamma$.
\end{longlist}
\end{proposition}

In the next theorems we assume that~\ref{hypoHgama} holds true
for some $\gamma$ in
$(0,2]$, and we study the self-similar limit of the rescaled solution
$w$ for $\mu^*={\mu}(\gamma)$.
We will see that the nondegeneracy of the limit will depend on the
shape of the spectral function $\mu$.
%
\begin{theorem}\label{teo2}
Let (\ref{H0-1}) be in force.
Assume that~\ref{hypoHgama} holds true for some $\gamma$ in $(0,2]$
and that ${\mu}(\delta)<{\mu}(\gamma)<+\infty$, for some $\delta
>\gamma$.
Then, there is a probability measure $\rho_{\infty,\gamma}$
such that:
\begin{longlist}[(iii)]
\item[(i)] The characteristic function of $\rho_{\infty,\gamma}$ is a
self-similar solution;
that is,
$\tphi_{\infty,\gamma}(\xi):=\int_\RE e^{i\xi v}\rho_{\infty
,\gamma}(dv)$
is a solution of (\ref{stazselfintegral}) for $\mu^*={\mu}(\gamma)$ and
%
\begin{equation}\label{lim-teo2}
\lim_{t \to+\infty} \phi\bigl(t,e^{-t\mu(\gamma)}\xi\bigr)
=\tphi_{\infty,\gamma}(\xi)
\end{equation}
for every $\xi\in\RE$.
Moreover,
\[
\tphi_{\infty,\gamma}(\xi)=\int_{\RE^+} \hat g_\gamma(\xi
z^{
{1}/{\gamma}}) \zeta_{\infty,\gamma}(dz),
\]
where
$\zeta_{\infty,\gamma}$ is given in \textup{(i)} of
Proposition~\ref{Pfixpoint}, and
$\hat g_\gamma$ is defined in (\ref{chaSta}).

\item[(ii)] If $\gamma\not=1,2$ or if $\gamma=1$ and
\textup{(b)} of~\ref{hypoH1} holds, then
$\rho_{\infty,\gamma}$ is a $\gamma$-stable distribution if and only
if
$\sum_{i=1}^N A_i^\gamma=1$ almost surely. Moreover\break $\int_\RE|v|^{p}
\rho_{\infty,\gamma}(dv)<+\infty$
if and only if $p<\gamma$.
\item[(iii)] If $\gamma=1$ and \textup{(a)} of~\ref{hypoH1} holds, then
$\rho_{\infty,\gamma}=\delta_{m_0}$ if and only
if $\sum_{i=1}^N A_i=1$ almost surely. Moreover, if $P\{\sum_{i=1}^N
A_i=1 \}<1$, then\break
$\int_\RE|v|^{p} \rho_{\infty,1}(dv)<+\infty$ for $p>1$
if and only if $p<q_{1}^*$ [where $q_1^*$ is defined in \textup{(ii)}
of Proposition~\ref{Pfixpoint}].

\item[(iv)] If $\gamma=2$, then $\rho_{\infty,2}$ is a Gaussian distribution
if and only
if $\sum_{i=1}^N A_i^2=1$ almost surely. Moreover, if $P\{\sum_{i=1}^N
A_i^2=1 \}<1$, then
$\int_\RE|v|^{p} \rho_{\infty,2}(dv)<+\infty$ for $p>2$
if and only if $p<q_{2}^*$ [where $q_2^*$ is defined in \textup{(ii)}
of Proposition~\ref{Pfixpoint}].
\end{longlist}
\end{theorem}

The previous result can be rephrased in terms of random variables as follows.
Let $V_t$ be a random variable whose characteristic function is the
unique solution $\phi(t,\xi)$ to problem
(\ref{eqboltzivp}). Such a $V_t$ is given explicitly in Section \ref
{Srandtreeandprobrep} (see\vspace*{1pt} Proposition~\ref{Propprobint}).
Then, Theorem~\ref{teo2} states the convergence in distribution of
$e^{-\mu(\gamma)t}V_t$ to a random variable $ V_\infty= Z_{\infty
,\gamma}^{{1}/{\gamma}} S_\gamma$,
where $Z_{\infty,\gamma}$ and $S_\gamma$ are stochastically
independent, $Z_{\infty,\gamma}$ is the unique
solution of (\ref{eqvRandomvariables}) with $\E[Z_{\infty
,\gamma}]=1$ and $ S_\gamma$ is
either a stable random variable or the constant~$m_0$.

The following theorem considers the cases in which the rescaling
$e^{-{\mu}(\gamma)t}$
provides a degenerate limiting solution, that is $V_\infty=0$.
%
\begin{theorem}\label{teo1}
Let (\ref{H0-1}) be in force. Assume that~\ref{hypoHgama} holds
true for some $\gamma$ in $(0,2]$.
If ${\mu}(\delta)<{\mu}(\gamma)<+\infty$, for some $0<\delta
<\gamma$, then
\[
\lim_{t \to+\infty} \phi\bigl(t,e^{-t\mu(\gamma)}\xi\bigr)=1 \qquad(\xi
\in\RE).
\]
\end{theorem}
%

\subsection{Comparison with previous results}\label{ssComparison}

In \cite{CeGaBo} the Chauchy problem (\ref{eqboltzivp})--(\ref{eq2})
is studied under the hypothesis that
$A_1,\ldots,A_N$ are exchangeable random variables with finite moments
of any order. The convergence of the
rescaled solution to the self-similar one is obtained under the
same hypotheses on the spectral function $\mu$ assumed in Theorem
\ref{teo2}, in two different situations:

- when $F_0$ is a symmetric distribution function and
%
\begin{equation}\label{condGAMBA}
1-\phi_0(\xi)=|\xi|^\gamma+ O(|\xi|^{\gamma+\eps}) \qquad(\xi
\to0)
\end{equation}
for some $\gamma\leq2$ and $\eps>0$;

- when $F_0$ is supported by $\RE^+$ and 
$L_0(\xi):=\int_{\RE^+}e^{-\xi v}\,dF_0(v)$
satisfies
%
\begin{equation}\label{condGAMBA2}
1-L_0(\xi)=\xi^\gamma+ O(|\xi|^{\gamma+\eps})\qquad(\xi\to0^+)
\end{equation}
for some $\gamma\leq1$ and $\eps>0$.

Some results on moments of the self-similar solutions
are proved under the stronger assumption that the distributions
of the $A_i$'s have compact support.

The probabilistic approach
enables us to weaken the hypotheses both on the $A_i$'s and on the
initial condition $\phi_0$.
In particular, we require neither the symmetry [except for \textup{(b)} in
assumption~\ref{hypoH1}]
nor the positiveness of the initial data. Moreover,
\ref{hypoHgama} is weaker than (\ref{condGAMBA}) and (\ref{condGAMBA2}).
In fact, if $F_0$ is symmetric, then it satisfies (\ref{stabledomain})
for $0<\gamma<2$
if and only if
\[
1-\phi_0(\xi)= k_0 |\xi|^\gamma
\bigl(1+o(1)\bigr)
\]
as $|\xi|\to0$, and $\s_0^2<+\infty$
if and only if $1-\phi_0(\xi)= \frac{\s^2_0}{2} |\xi|^2(1+o(1)) $.
See Theorem 1.3 of \cite{Ibragimov1985}.
On the other hand, if $F_0$ is supported by $\RE^+$ and its Laplace transfom
satisfies (\ref{condGAMBA2}) for $\gamma<1$, by Theorem 4 in Section
XII.5 of \cite{FellerVolII}, it follows that
(\ref{stabledomain}) holds true.
Finally, if (\ref{condGAMBA2}) holds for $\gamma=1$,
it follows immediately that $\int_{\RE^+}v \,dF_0(v)<+\infty$.

\subsection{Rates of convergence}
In the following
we show that, under some additional hypotheses, the convergence stated
in Theorem~\ref{teo2} takes place
at an exponential rate in suitable metrics.

Recall that the Wasserstein distance of order $\delta>0$
between two random variables $X$ and $Y$, or equivalently between their
probability distributions,
is defined by
%
\begin{equation}
\label{eqwasserstein}
\W_\delta(X,Y):=\inf_{(X',Y')}(\E|X'-Y'|^\delta)^{{1}/{\max
(\delta
,1)}} .
\end{equation}
The infimum is taken over all pairs $(X',Y')$ of real random variables
whose marginal distributions are the same as those of $X$ and $Y$, respectively.
In general, the infimum in (\ref{eqwasserstein}) may be infinite;
a sufficient (but not necessary) condition for finite distance is
that both $\E|X|^\delta<+\infty$ and $\E|Y|^\delta<+\infty$. For more
information on Wasserstein distances see, for example, \cite{Rachev1991}.

The next theorem is the natural generalization of Theorem 5 in
\cite{BLM2008}.
%
\begin{theorem}
\label{PropW-2}
Let (\ref{H0-1}) be in force. Assume that~\ref{hypoHgama} holds
true for some
$\gamma$ in $(0,2)$ and that ${\mu}(\delta)<{\mu}(\gamma)$, for some
$\gamma< \delta$
with $1\leq\gamma<\delta\leq2$
or $\gamma<\delta\leq1$. Let $V_t$ and $V_\infty$ be as above.
Then
%
\begin{equation}
\label{mainbound}
\W_{\delta}\bigl(e^{-\mu(\gamma)t}V_t, V_\infty\bigr)^{\max(\delta,1)}
\leq c \W
_{\delta}
(X_0,V_\infty)^{\max(\delta,1)}
e^{-t \delta[{\mu}(\gamma)-{\mu}(\delta)]}
\end{equation}
with $c=1$ if $\delta\leq1$ and $c=2$ otherwise.\vadjust{\goodbreak}
\end{theorem}

Clearly, (\ref{mainbound}) is meaningful only if $\W_{\delta
}(X_0,V_\infty)<+\infty$.
When $\gamma=1$ and (a) of~\ref{hypoH1} holds, it follows that
$\E| V_\infty|^\delta<+\infty$ by (iii) of Theorem~\ref{teo2},
since it is easy to see that $\delta<q^*_1$. Hence, in this case, $\W
_{\delta}(X_0,V_\infty)<+\infty$
whenever $\E|X_0|^{\delta}<+\infty$.
When $\gamma\not=1$ or when $\gamma=1$ and (b) of~\ref{hypoH1} holds, the requirement
$\W_{\delta}(X_0,V_\infty) < +\infty$ is nontrivial
since, by Theorem~\ref{teo2}, one has $\E[|V_\infty|^\delta]=+
\infty$.
The following lemma provides a sufficient criterion tailored to the
situation at hand.
%
\begin{lemma}
\label{lemma-tail2}
Assume, in addition to the hypotheses of Theorem~\ref{teo2},
that $\delta<2\gamma$
and that $F_0$
satisfies~\ref{hypoHgama} in the more restrictive sense that
there exists a constant $K>0$ and some $0<\varepsilon<1$ with
%
\begin{eqnarray}
\label{eqFtail}
|1-c_0^+x^{-\gamma} - F_0(x)| & < & K x^{-(\gamma+\varepsilon)}
\qquad\mbox
{for $x>0$}, \\
\label{eqFtail2}
|F_0(x) - c_0^-(-x)^{-\gamma}| & < & K (-x)^{-(\gamma+\varepsilon)}
\qquad\mbox{for $x<0$}.
\end{eqnarray}
Provided that $\delta<\gamma/(1-\eps)$, it follows $\W_{\delta
}(X_0,V_\infty)<+\infty$.
\end{lemma}

We have not been able to prove Theorem~\ref{PropW-2} for $\gamma=2$. On
the other hand
we are able to give the speed of convergence for every $\gamma$ in
$(0,2]$ with respect
to the Zolotarev metrics $\ZZ_s$.
The metric $\ZZ_s$ is defined, for $s=m+\alpha$, $m$ being a
nonnegative integer and $0<\alpha\leq1$,
by
%
\begin{equation}
\label{eqzolo}
\ZZ_s(X,Y):=\sup\{ \E[ f(X)-f(Y)] \dvtx f \in\CF_s \},
\end{equation}
where
$\CF_s$ is the set of real-valued functions on $\RE$ which at all
points have
the $m$th derivatives such that
$|f^{(m)}(x)-f^{(m)}(y)|\leq|x-y|^\alpha$. For more information see,
for example, \cite{Zolometric}.

In general the finiteness condition $\ZZ_s(X,Y)$ is not easy to check.
It turns out that if
$\int x^r(dF_X(x)-dF_Y(x)) =0$ for any integer $r \leq m$ and
$\int|x|^s |dF_X(x)-dF_Y(x)| <+\infty$, where $F_X$ and $F_Y$
are the distribution functions of $X$ and~$Y$,
then $\ZZ_s(X,Y)<+\infty$.
See Theorem 1.5.7 in \cite{Zolometric}.
The estimates proved in the next
theorem are interesting in particular for the case $\gamma=2$, for which
the above sufficient conditions for the finiteness of $\ZZ
_s(X_0,V_\infty)$ are easily verified.
%
\begin{theorem}
\label{PropZolo}
Let (\ref{H0-1}) be in force. Assume that~\ref{hypoHgama} holds
true for some
$\gamma$ in $(0,2]$ and that ${\mu}(\delta)<{\mu}(\gamma)$, for some
$\gamma< \delta$.
Let $V_t$ and $V_\infty$ be as above.
Then
%
\begin{equation}
\label{mainbound2}
\ZZ_{\delta}\bigl(e^{-\mu(\gamma)t}V_t, V_\infty\bigr) \leq\ZZ_{\delta}
(X_0,V_\infty) e^{-t \delta[{\mu}(\gamma)-{\mu}(\delta)]}.
\end{equation}
In particular, if $\gamma=2$, $\delta\leq3$ and $\E|X_0|^\delta
<+\infty
$, then
%
\begin{equation}
\label{mainbound3}
\ZZ_{\delta}\bigl(e^{-\mu(\gamma)t}V_t, V_\infty\bigr) \leq
c e^{-t \delta[{\mu}(\gamma)-{\mu}(\delta)]},
\end{equation}
where
\[
c := \frac{1}{\Gamma(1+\delta)}(\E|X_0|^\delta+E|V_\infty
|^\delta
)<+\infty.
\]
\end{theorem}
%

\section{Marked recursive $N$-ary random trees and probabilistic
interpretation of the solutions}\label{Srandtreeandprobrep}

The notion of $N$-ary random trees will be used to describe, in a
probabilistic way,
the solution of (\ref{eqboltzivp}). This approach is a generalization
of the probabilistic representation presented in \cite{BLM2008}, where
binary trees were considered
in order to describe the solution when $\hat Q$ is a bilinear smoothing
transformation.

\subsection{Random $N$-ary recursive trees}\label{Snary}
Recall that a rooted tree is said to be a planar tree
when successors of the root and recursively the successors of
each node are equipped with a left-to-right order.
For any integer number $N \geq2$ an $N$-ary tree is
a (planar and rooted) tree
where each node is either a leaf (i.e., it
has no successor) or it has $N$ successors.
We define the size of the $N$-ary tree~$t$, in symbol $|t|$, by the
number of internal nodes. Any $N$-ary tree with $Nk+1$ nodes
has size $k$ and possesses $f_k:=(N-1)k+1$ leaves.
We now describe a (natural) tree evolution process which gives rise
to the so-called ``random $N$-ary recursive tree.''
The evolution process starts with~$T_0$, an empty tree,
that is, with just an external node (the root). The first step in the
growth process is
to replace this external node by an internal one with~$N$ successors that
are leaves; in this way we get $T_1$. Then with probability $1/N$ (i.e.,
the number of leaves)
one of these $N$ leaves is selected and again replaced by an internal
node with~$N$ successors. In this way one continues.
At every time $k$, $T_k$ is an $N$-ary tree with $k$ internal nodes.

A very important issue is that $N$-ary trees
have a recursive structure. More precisely we can use the following
recursive definition of $N$-ary trees: an $N$-ary tree $t$ is either
just an external node or an internal node
with~$N$ subtrees that are again $N$-ary trees.
We shall denote these subtrees by $t^{(1)},\ldots,t^{(N)}$.

Recall also that every $N$-ary tree
can be seen as a subset of
\[
\U:=\{\varnothing\}\cup\biggl[\bigcup_{k \geq1 } \{1,2,\ldots,N\}^k\biggr].
\]
As usual $\varnothing$ is the root, and if
$v=(v_1,\ldots,v_k)$ ($v_i \in\{1,\ldots,N\})$ is a node of an $N$-ary
tree then
the length of $v$ is $|v|:=k$, moreover $(v,v_{k+1}):=(v_1,\ldots
,v_k,v_{k+1})$, and $(v,\varnothing):=v$.
For every $1 \leq i\leq k$, set
$v|i:=(v_1,\ldots,v_i)$
and $v|0=\varnothing$.
Finally, given an $N$-ary tree $t$ we shall denote
by $\CL(t)$
the set of the leaves of $t$.
For more details on $N$-ary recursive trees see, for instance,~\cite{drmota}.

For every integer $k \geq1$, set
\[
\CI_k:= \Biggl\{ \underline{i}=(i_1,\ldots,i_N) \in\N_0^N\dvtx \sum_{j=1}^N
i_j=k-1 \Biggr\},
\]
where $\N_0=0 \cup\N$, and denote by $\T_{k,N}$ the set of all $N$-ary
trees with size~$k$. Notice that $\underline{i}\in\CI_1$ if, and only
if, $i_1=\cdots= i_N=0$.\vadjust{\goodbreak}
Finally, for every~$\underline{i}$ in~$\CI_k$, set
\[
C^k_{\underline{i}}=\bigl\{ t \in\T_{k,N} \dvtx \bigl|t^{(j)}\bigr|=i_j,    j=1,\ldots
,N\bigr\}.
\]

The following proposition states some properties of random $N$-ary
recursive trees.\vspace*{-2pt}
%
\begin{proposition}\label{comblemma} Let $(T_k)_{k \geq1}$ be a
sequence of random $N$-ary recursive trees.
For every $k \geq1$, every~$\underline{i}$ in $\CI_k$ and every $t$ in
$\T_{k,N}$,
%
\begin{eqnarray}\label{fact-tree}
&&
P\bigl\{ T^{(1)}_k=t^{(1)},\ldots, T^{(N)}_k=t^{(N)}| T_k \in
C^k_{\underline{i}}\bigr\} \nonumber\\[-8pt]\\[-8pt]
&&\qquad =\prod_{j=1}^N P\bigl\{ T_{|t^{(j)}|} =t^{(j)}\bigr\}
\J\bigl\{\bigl|t^{(j)}\bigr|=i_j\bigr\}
\nonumber
\end{eqnarray}
and
%
\begin{equation}\label{fact-tree2}
P\{T_k \in C^k_{\underline{i}}\}=p_k(\underline{i}),
\end{equation}
where for $k\geq1$,
%
\begin{equation}\label{espr-p}
p_k(\underline{i}):= 
\pmatrix{k-1 \cr i_1,\ldots,i_N}
\frac{\prod_{l=1}^N\prod_{m=0}^{i_l-1}f_m}{\prod_{r=0}^{k-1}f_r}
\end{equation}
with the convention that ${\prod_{r=0}^{-1}f_r}=1$.
Finally,
$(|T_{n}^{(1)}|/n,\ldots,|T_{n}^{(N)}|/n)$ converges almost surely
to a vector $(U_1,\ldots,U_N)$
with Dirichlet distribution of parameters $(1/(N-1),\ldots,1/(N-1))$,
for $n \to+\infty$.\vspace*{-2pt}
\end{proposition}

\subsection{Wild series and probabilistic representation of the solutions}
To start with we will write the Wild series expansion of $\phi(\cdot,t)$.
This kind of expansion can be easily derived using
a general result contained in \cite{Kielek88}.
For every $t \geq0$ and $k \in\N_0$, set
\[
b_k:=\frac{\prod_{i=0}^{k-1} f_i }{(N-1)^kk!}
\]
and
%
\begin{equation}\label{defnegativebinomial}
\zeta(t,k):=b_k e^{-t}\bigl(1-e^{-(N-1)t} \bigr)^k.\vspace*{-2pt}
\end{equation}

\begin{remark}\label{remark1}
Note that $\zeta(t,\cdot)$ is the probability density of a negative binomial
random variable of parameters $(1/(N-1),e^{-(N-1)t})$. Indeed, since $f_i=
(N-1)i+1$ for $i=0,1, \ldots$
\[
b_k=\frac{(1/(N-1))_k}{k!},
\]
where for every nonnegative real number $r$, and every nonnegative
integer~$n$
\[
(r)_n=\prod_{i=0}^{n-1}(r+i)=\frac{\Gamma(r+n)}{\Gamma(r)}
\]
and $(r)_0=1$.\vadjust{\goodbreak}
\end{remark}

Using Remark~\ref{remark1} above and Theorem 1.1 in \cite{Kielek88}
it is a simple matter to deduce that
the unique global solution to (\ref{eqboltzivp})
is given by
\[
\phi(t,\xi)=\sum_{k \geq0} \zeta(t,k) q_k(\xi),
\]
where
$(q_k)_k$ is a sequence of characteristic functions
recursively defined by setting $q_0(\xi)=\phi_0(\xi)$ and, for
$k\geq1$,
\[
q_k(\xi)=\sum_{ \underline{i} \in\CI_k} p_k(\underline{i})
\hat
{Q}(q_{i_1},\ldots,q_{i_N})(\xi),
\]
where $p_k$ is defined in (\ref{espr-p}).
This representation is the generalization of the Wild series, which is
obtained, when $N=2$, in \cite{Wild1951}.
It is easy to see that $\phi(t,\cdot)$ is a characteristic function.

The Wild series expansion suggests a probabilistic interpretation for
the solutions
as sums of random variables indexed by $N$-ary recursive random trees.
On a sufficiently large probability space $(\Omega,\CF,P)$,
let the following be given:
\begin{itemize}
\item a family $(X_v)_{v \in\U}$ of independent random variables
with common distribution function $F_0$;
\item a family $( A_1(v),A_2{(v)},\ldots, A_N{(v)} )_{v \in
\U}$
of independent positive random
vectors with the same distribution of $(A_1,\ldots,A_N)$;
\item a sequence of $N$-ary recursive random trees $(T_n)_{n\in\setN}$;
\item a stochastic process $(\nu_t)_{t\geq0}$ with values in $\N_0$
such that
$P\{\nu_t=k\}=\zeta(t,k)$ for every integer $k \geq0$, where $\zeta
(t,k)$ is defined in (\ref{defnegativebinomial}).
\end{itemize}
Write $A{(v)}= (A_1(v),A_2{(v)},\ldots, A_N{(v)})$ and
assume further that
\[
(A{(v)})_{v \in\U},\qquad (T_n)_{n \geq1},\qquad (X_v)_{v \in\U}
\quad\mbox{and}\quad
(\nu_t)_{ t>0 }
\]
are stochastically independent.

For each node $v=(v_1,\ldots,v_k)$ in $\U$, set
\[
\be(v):=\prod_{i=0}^{|v|-1} A_{v_{i+1}}(v|i)
\]
and $\be(\varnothing)=1$.
Now recall that $\CL(T_n)$ is the set of leaves of $T_n$, and define
\[
W_0:=X_\varnothing
\]
and, for every $n \geq1$,
\[
W_{n}:=\sum_{v \in\CL(T_{n})} \be(v) X_v.
\]
%
\begin{proposition}\label{Propprobint}
Equation (\ref{eqboltzivp}) has a unique solution $\phi(t,\cdot)$, which
coincides with the characteristic function of
$ V_t:=W_{\nu_t}$.
\end{proposition}

Let us conclude this section rewriting $W_n$ in an alternative form.
In the following we will use both forms, according to our convenience.
For each $n \geq1$ we shall denote by
\[
\{ \beta_{1,n},\ldots, \beta_{f_n,n} \}
\]
the weights associated to the leaves of $T_n$, that is, if
\[
\CL(T_{n})=\{L_{1,n},\ldots, L_{f_n,n} \}
\]
(in left-to-right order) $\beta_{i,n}=\be(L_{i,n})$. Hence we can rewrite
$W_n$ as
\[
W_n=\sum_{j=1}^{f_n} \beta_{j,n} X_{j,n},
\]
where
$
X_{j,n}:=X_{L_{j,n}}$.

\section{Some limit theorems for sums of random variables indexed by
$N$-ary recursive trees}\label{Smartingaleetal}
Let us sketch our approach to the study of the asymptotic behavior of
$\phi(t,e^{-\mu(\gamma)t}\xi)$.
From the probabilistic interpretation given in Proposition \ref
{Propprobint}, we obtain that
$\phi(t,e^{-\mu(\gamma)t}\xi)$ is the characteristic function
of the rescaled random variable $e^{-\mu(\gamma)t} V_t=e^{-\mu
(\gamma
)t} W_{\nu_t}$.
Hence, we look for a~positive function $n \mapsto m_{n}(\gamma)$
such that
%
\begin{equation}\label{vettore}
(N_t(\gamma),\tilde W_{\nu_t}):=\Biggl(e^{-\mu(\gamma)t}
m_{\nu_t}(\gamma)^{{1}/{\gamma}}, \sum_{j=1}^{f_{\nu_t}}
\frac{\beta_{j,\nu_t}}{m_{\nu_t}(\gamma)^{{1}/{\gamma}}}
X_{j,\nu_t}
\Biggr)
\end{equation}
converges\vspace*{1pt} weakly as $t \to+\infty$, in order to obtain the
convergence of $ e^{-\mu(\gamma)t}V_t=N_t(\gamma) \tilde W_{\nu_t} $.
This will be done in several steps. First of all we will study, for
suitable $m_{n}(\gamma)$'s, the weak limit of
%
\begin{equation}\label{defWnbis}
\tilde W_n:=\frac{W_n}{m_{n}(\gamma)^{{1}/{\gamma}}} = \frac
{1}{m_{n}(\gamma)^{{1}/{\gamma}}}
\sum_{j=1}^{f_{n}} \beta_{j,n} X_{j,n},
\end{equation}
which is a sum of random variables from a triangular array. Notice that
a~direct application of a central limit theorem is not possible,
since the weights $m_{n}(\gamma)^{-{1}/{\gamma}}\beta_{j,n}$
are not independent. However, one can apply a central limit theorem to
the conditional
law of $\tilde W_n$, given the array of weights $\beta_{j,n}$ and
$(T_n)_{n \geq1}$.
To this end, we shall prove that under suitable assumptions, if
%
\begin{equation}\label{defmngamma}
m_n(\gamma):=\prod_{k=0}^{n-1} \biggl(1 +\frac{\CS(\gamma
)}{f_k}\biggr),
\end{equation}
$\QQ(\gamma)$ being defined in (\ref{eqthes}),
then
\[
\tilde M_n(\gamma):= \frac{1}{m_n(\gamma)} \sum_{j=1}^{f_n} \beta
_{j,n}^\gamma
\]
converges a.s. to a limit $\tilde M_{\infty}(\gamma)$ and
that $\max_{j=1,\ldots,f_n} \beta_{j,n}m_n(\gamma)^{-
{1}/{\gamma}}$
converges to
zero in probability as $n \to+\infty$.
As a consequence we will find that the weak limit of $\tilde W_n$
is a scale mixture of $\gamma$-stable laws, where the scale mixing measure
is the law of $\tilde M_{\infty}(\gamma)^{{1}/{\gamma}}$. From the
asymptotic results on
$\tilde W_n$ we will easily deduce the weak convergence of
the random vector (\ref{vettore}) for $t \to+\infty$.

%
\subsection{The martingale of weights}\label{Smartingweights}
Let $\gamma$ be a given positive real number such that
$\E[\sum_{i=1}^N A_i^\gamma]<+\infty$.
For every integer number $n \geq1$, set
%
\begin{equation}
\label{eqthem}
M_n{(\gamma)} := \sum_{v \in\CL(T_n)} \be(v)^\gamma=\sum_{j=1}^{f_n}
\beta_{j,n}^\gamma.
\end{equation}
Note that
\[
\tilde M_n (\gamma)= \frac{M_n (\gamma)}{m_n(\gamma)}
\]
and, if $\CS(\gamma)=0$, then $\tilde M_n (\gamma)= M_n (\gamma)$.

The following proposition generalizes Lemma 2 in \cite{BLM2008}.
%
\begin{proposition}\label{Lemma2}
For every $\gamma>0$ such that $\E[\sum_{i=1}^N A_i^\gamma]
<+\infty$,
one~has
\[
\E[M_n(\gamma)]= m_n(\gamma)=\frac{(({\CS(\gamma
)+1})/({N-1})
)_n}{({1}/({N-1}))_n}
\]
and, as $n \to+\infty$,
%
\begin{equation}\label{asintm}
m_n(\gamma)=n^{{\QQ(\gamma)}/({N-1})} \frac{\Gamma(
{1}/({N-1}))}{\Gamma(({\QQ(\gamma)+1})/({N-1}))}
\biggl(1+O\biggl(\frac{1}{n}\biggr)\biggr).
\end{equation}
Moreover, $\tilde M_n(\gamma)$ is a positive martingale with respect to
the filtration
\[
\CG_n=\s\bigl((A{(v)})_{v \in\CI(T_{n})},T_1,\ldots,T_{n}\bigr),
\]
where\vspace*{1pt} $\CI(T_{n})$ denotes the set of the internal nodes of $T_n$,
and \mbox{$\E[\tilde M_n(\gamma)]=1$}. Hence, $\tilde M_n(\gamma)$ converges
almost surely to
a random variable $\tilde M_\infty(\gamma)$ with $\E[\tilde M_\infty
(\gamma)] \leq1$.
\end{proposition}

For every $\gamma>0$, set
\[
\b_{(n)}^{(\gamma)}:=\max_{v \in\CL(T_n)} \frac{\be
(v)}{m_n(\gamma
)^{{1}/{\gamma}}}
=\max_{j=1,\ldots,f_n} \frac{\beta_{j,n}}{m_n(\gamma)^{
{1}/{\gamma}}}
\]
and recall that
${\mu}(\gamma)=\QQ(\gamma)/\gamma$.
%
\begin{proposition}\label{Lemma4}
If for some $\delta>0$ and $\gamma>0$ one has
${\mu}(\delta) < {\mu}(\gamma)<+\infty$, then
$\b_{(n)}^{(\gamma)}$ converges in probability to $0$.
Moreover, if in addition $\delta<\gamma$ one has
that $\tilde M_n(\gamma)$ converges almost surely to $0$, that is,
$\tilde M_\infty(\gamma)=0$.
While, if $\gamma< \delta$, $\tilde M_n(\gamma)$ converges in $L^1$ to
$\tilde M_\infty(\gamma)$
and hence $\E[\tilde M_\infty(\gamma)]=1$.
\end{proposition}
%
\begin{proposition}\label{propeqlimiteMinfty}
Assume that $\E[\sum_{i=1}^N A_i^\gamma]<+\infty$. Let
$\tilde M_\infty(\gamma)$ be the same random variable defined in
Proposition~\ref{Lemma2},
and denote its characteristic function by $\psi_{\infty,\gamma}$.
Then
$\psi_{\infty,\gamma}$
satisfies the following integral equation:
%
\begin{equation}\label{integraleqpsi}
\psi_{\infty,\gamma}(\xi)=\E\Biggl[ \prod_{i=1}^N \psi_{\infty
,\gamma
}\bigl(A_i^\gamma U_i^{{\QQ(\gamma)}/{(N-1)}}\xi\bigr) \Biggr] \qquad(\xi
\in
\RE),
\end{equation}
where $U=(U_1,\ldots,U_N)$ has
Dirichlet distribution of parameters $(1/(N-1),\allowbreak\ldots,1/(N-1))$
and $(A_1,\ldots,A_N)$ and $U$ are stochastically independent.
\end{proposition}

Note that (\ref{integraleqpsi}) is equivalent to
%
\begin{equation}\label{integraleqpsi2}
M\stackrel{\CL}=\sum_{i=1}^N A_i^\gamma U_i^{{\QQ(\gamma
)}/{(N-1)}} M_i,
\end{equation}
where $(M,M_1,\ldots,M_n)$ are stochastically independent random variables
with the same law of $\tilde M_\infty(\gamma)$, and
$(M,M_1,\ldots,M_n)$, $U$ and $(A_1,\ldots,A_N)$
are stochastically independent.

\subsection{\texorpdfstring{Convergence of $\tilde W_n$ and of $(N_t(\gamma),\tilde W_{\nu_t})$}
{Convergence of W n and of (N t(gamma), W nu t)}}\label{jointconv}

Now we study the limiting distribution of $\tilde W_n$
defined by (\ref{defWnbis}).
%
\begin{proposition}\label{proplimit1}
Let (\ref{H0-1}) be in force. Let $\gamma$ belong to $(0,2]$, and
assume that there exists $\delta>0$ such that ${\mu}(\delta)<{\mu
}(\gamma
)<+\infty$.
Assume that condition~\ref{hypoHgama} holds true; then
%
\begin{equation}
\label{eqchar02}
\lim_{n \to+\infty} \E[ e^{ i \xi\tilde W_{n}}]
=\E[ \hat g_\gamma( \tilde M_{\infty}(\gamma)^{{1}/{\gamma}}
\xi) ]
\end{equation}
for every $\xi\in\RE$, where $\tilde M_\infty(\gamma)$ is the same
random variable
defined in Proposition~\ref{Lemma2} and $\hat g_\gamma$ is
defined in (\ref{chaSta}).
\end{proposition}

At this stage, recall that $N_t(\gamma)=e^{-\mu(\gamma)t}m_{\nu
_t}(\gamma)^{1/\gamma}$.
%
\begin{proposition}\label{Propjointconv}
Under the assumptions of Proposition~\ref{proplimit1},
\begin{eqnarray*}
&&
\lim_{t \to+\infty} \E\bigl[e^{i\xi_1 N_t(\gamma) +i\xi_2 \tilde
W_{\nu_t}}\bigr]\\
&&\qquad=\E\bigl[e^{i\xi_1 c_\gamma\GA^{{{\mu}(\gamma)}/{(N-1)}} } \bigr]
\E[ \hat g_\gamma( \tilde M_{\infty}(\gamma)^{{1}/{\gamma}}
\xi_2) ],
\qquad(\xi_1,\xi_2) \in\RE^2,
\end{eqnarray*}
where $\GA$ has gamma distribution with shape parameter $1/(N-1)$ and
scale parameter $1$,
%
\begin{equation}\label{defcgamma}
c_\gamma:=\biggl( \frac{\Gamma({1}/({N-1}))}{\Gamma(({\QQ
(\gamma
)+1})/({N-1}))} \biggr)^{{1}/{\gamma}},
\end{equation}
$\tilde M_\infty(\gamma)$ is the same random variable
defined in Proposition~\ref{Lemma2}, $\tilde M_\infty(\gamma)$
and~$\GA$ are stochastically independent and
$\hat g_\gamma$ is defined in (\ref{chaSta}).
As a consequence,
%
\begin{equation}\label{limitw}\hspace*{29pt}
\lim_{t \to+\infty} \phi\bigl(t,e^{-\mu(\gamma)t}\xi\bigr)=
\E\bigl[ \hat g_\gamma\bigl( c_\gamma\GA^{{{\mu}(\gamma)}/{(N-1)}}
\tilde
M_{\infty}(\gamma)^{{1}/{\gamma}} \xi\bigr) \bigr]
\qquad(\xi\in\RE).
\end{equation}
\end{proposition}

The result in (\ref{limitw}) is the core of
Theorems~\ref{teo2} and~\ref{teo1} presented in Section~\ref{SecMain1}.
The further properties of the limiting distribution
are proved in Section~\ref{proofsMain}.\vspace*{-2pt}

\section{Proofs}\label{Sproofs}\vspace*{-2pt}

\subsection{\texorpdfstring{Proofs of Section \protect\ref{Srandtreeandprobrep}}{Proofs of Section 3}}\vspace*{-2pt}

\mbox{}

\begin{pf*}{Proof of Proposition~\ref{comblemma}}
Let us first prove (\ref{fact-tree2}).
Recall that $\CI_1=\{(0,\ldots,0) \}$ and
for $\underline{i}=(0,\ldots,0)$
\[
P\{T_1 \in C^1_{\underline{i}}\}=P\bigl\{\bigl|T_1^{(1)}\bigr|=0,\ldots,
\bigl|T_1^{(N)}\bigr|=0\bigr\}=1
=p_1(\underline{i}).
\]
For every $k\geq1$ and $j=1,\ldots,N$,
\begin{eqnarray*}
&&
P\bigl\{\bigl|T_{k+1}^{(l)}\bigr| = \bigl|T_{k}^{(l)}\bigr|,
l\neq j,
\bigl|T_{k+1}^{(j)}\bigr|=\bigl|T_{k}^{(j)}\bigr|+1
\big|\bigl|T_{k}^{(l)}\bigr|, l=1,\ldots,N \bigr\} \\
&&\qquad =
\frac{f_{|T_{k}^{(j)}|}}{f_k}.
\end{eqnarray*}
This means that the problem of evaluating probability (\ref
{fact-tree2}) can be reduced to a P\'olya urn scheme, where
one starts with~$N$ balls of $N$ different colors, and at each step a
ball is randomly drawn from the urn and replaced
with~$N$ balls of the same color. Hence, for every $k\geq2$ and $\underline
{i}=(i_1, \ldots, i_N) \in\CI_{k}$
\[
P\bigl\{\bigl|T_{k}^{(l)}\bigr|= i_l, l=1,N\bigr\}=\frac{(k-1)!}{\prod
_{l=1}^{N}i_l!}
\frac{\prod_{l=1}^{N}\prod_{m=0}^{i_l-1}f_m}{\prod
_{r=0}^{k-1}f_r}=p_k(\underline{i}),
\]
which is (\ref{fact-tree2}).

Let us prove (\ref{fact-tree}) by induction. For $k=1$ equality (\ref
{fact-tree}) is trivially true. Let us suppose (\ref{fact-tree}) holds
for $k$.
Let $t \in\T_{k+1,N}$ and $\underline{i}=(i_1, \ldots,
i_N)=(|t^{(1)}|,\ldots,|t^{(N)}|)\in\CI_{k+1}$, then
\begin{eqnarray*}
\hspace*{-4.5pt}&& P\bigl\{ T^{(1)}_{k+1}=t^{(1)},\ldots, T^{(N)}_{k+1}=t^{(N)}, T_{k+1} \in
C^{k+1}_{\underline{i}}\bigr\} \\
\hspace*{-4.5pt}&&\qquad= P\bigl\{ T^{(1)}_{k+1}=t^{(1)},\ldots, T^{(N)}_{k+1}=t^{(N)}\bigr\} \\
\hspace*{-4.5pt}&&\qquad= \hspace*{-0.7pt}\mathop{\sum_{j=1,\ldots,N}}_{j\dvtx i_j \geq1}\sum_{t^*_j\in A_{j,t}}
\hspace*{-1.8pt}P\bigl\{ T^{(1)}_{k+1}=t^{(1)},\ldots,
T^{(N)}_{k+1}=t^{(N)}|T^{(l)}_{k}=t^{(l)}   ,l\neq j,
T^{(j)}_{k}= t^*_{j}\bigr\} \\
\hspace*{-4.5pt}&&\qquad\quad\hspace*{60.2pt}{}\times P\bigl\{
T^{(l)}_{k}=t^{(l)}  , l\neq j,    T^{(j)}_{k}= t^*_{j}\bigr\},
\end{eqnarray*}
where
\[
A_{j,t}=\bigl\{t^*_j\in\T_{i_j-1,N} \dvtx\exists   v\in\CL(t^*_j)
\mbox{ such that }
t^*_j\cup\{(v,1), \ldots, (v,N)\}=t^{(j)} \bigr\}.
\]
By construction of a random $N$-ary tree, if $i_j \geq1$ and $t^*_j
\in A_{j,t}$,
%
\begin{equation}\label{pezzo1}
P\bigl\{ T^{(1)}_{k+1} =t^{(1)},\ldots,
T^{(N)}_{k+1}=t^{(N)}|T^{(l)}_{k}=t^{(l)}  , l\neq j,
T^{(j)}_{k}= t^*_{j}\bigr\}
=\frac{1}{f_k}
\end{equation}
and
%
\begin{equation}\label{pezzo2}
P\bigl\{T_{i_j}=t^{(j)}\bigr\}=\frac{1}{f_{i_j-1}} \sum_{t^*_j\in A_{j,t}}P\{
T_{i_j-1}=t^*_j\}.
\end{equation}
Furthermore, in view of the induction hypotheses and (\ref
{fact-tree2}), one gets
%
\begin{eqnarray}\label{pezzo3}
&&
P\bigl\{T^{(l)}_{k}=t^{(l)}  , l\neq j,    T^{(j)}_{k}=
t^*_{j}\bigr\}\nonumber\\[-8pt]\\[-8pt]
&&\qquad=\prod_{l=1}^N P\bigl\{T_{i_l}=t^{(l)}\bigr\}
\frac{P\{T_{i_j-1}=t^*_j\}}{P\{
T_{i_j}=t^{(j)}\}}
p_k(i_1,\ldots,i_{j}-1,\ldots,i_N).\nonumber
\end{eqnarray}
Hence, from (\ref{pezzo1}), (\ref{pezzo2}) and (\ref{pezzo3}), one obtains
\begin{eqnarray*}
&&
P\bigl\{ T^{(1)}_{k+1}=t^{(1)},\ldots, T^{(N)}_{k+1}=t^{(N)}\bigr\}\\
&&\qquad
=\prod_{l=1}^N P\bigl\{T_{i_l}=t^{(l)}\bigr\}\mathop{\sum_{j=1,\ldots,N}}_{j\dvtx i_j
\geq1}
\frac{p_k(i_1,\ldots,i_{j}-1,\ldots,i_N)}{f_k}
\sum_{t^*_j\in A_{j,t}}\frac{P\{T_{i_j-1}=t^*_j\}}{P\{
T_{i_j}=t^{(j)}\}
} \\
&&\qquad
= \prod_{l=1}^N
P\bigl\{T_{i_l}=t^{(l)}\bigr\}\mathop{\sum_{j=1,\ldots,N}}_{j\dvtx
i_j \geq1} \frac{f_{i_j-1}}{f_k}p_k(i_1,\ldots,i_{j}-1,\ldots,i_N) \\
&&\qquad= \prod_{l=1}^N P\bigl\{T_{i_l}=t^{(l)}\bigr\} p_{k+1}(i_1,\ldots,i_N),
\end{eqnarray*}
where the last equality is obtained by direct replacement of the
expression of $p_k(i_1,\ldots,i_{j}-1,\ldots,i_N)$.
Note that, using the P\'olya urn interpretation, $|T_k^{(l)}|$
represents the numbers of
balls of color $l$ drawn in the first $k-1$ steps.
Hence, using the results in \cite{blackwellMacQueen},
the almost sure convergence of $(|T_k^{(l)}|/(k-1)\dvtx\break l=1,\ldots,N)$
follows by the strong law of large numbers for exchangeable sequences.
\end{pf*}
\begin{pf*}{Proof of Proposition~\ref{Propprobint}}
We
need only to prove that
$q_n(\xi)=\E[e^{i\xi W_n}]$, for every $n\geq0$. This is clearly true
for $n=0$.
For $n\geq1$, write
\[
W_n=\sum_{j=1}^{N} A_j{(\varnothing)} \Biggl[ \sum_{v \in\CL(T_n^{(j)})}
\prod_{i=0}^{|v|-1} A^{(j)}_{v_{i+1}}(v|i) X^{(j)}_v \Biggr],
\]
where $A^{(j)}(v)=A((j,v))$, $ X^{(j)}_v=X_{(j,v)}$, and, by
convention, if $\CL(T_n)=\varnothing$
the term between square brackets is equal to $X_{(j)}$. Since
$(A^{(j)}(v), X^{(j)}_v)_{v\in\U}$, $j=1,\ldots,N$, are
independent, with the same distribution of $(A(v), X_v)_{v\in\U}$,
using (\ref{fact-tree}) and the
induction hypothesis, one proves that
%
\begin{equation}\label{fact2}
\E\bigl[e^{i\xi W_n}|A{(\varnothing)},\bigl|T_n^{(1)}\bigr|,\ldots,
\bigl|T_n^{(N)}\bigr|\bigr]=
\prod_{j=1}^N q_{|T_n^{(j)}|}(\xi A_j{(\varnothing)}).
\end{equation}
At this stage the conclusion follows easily by using (\ref{fact-tree2});
indeed,
\[
\E[e^{i\xi W_n}]=\E\Biggl[\prod_{j=1}^N q_{|T_n^{(j)}|}(\xi
A_j{(\varnothing)})\Biggr]
=\sum_{\underline{i} \in\CI_n} \E\Biggl[\prod_{j=1}^N q_{i_j}(\xi
A_j)\Biggr] p_{n}(\underline{i})
=q_n(\xi).\quad
\]
\upqed
\end{pf*}

\subsection{\texorpdfstring{Proofs of Section \protect\ref{Smartingaleetal}}{Proofs of Section 4}}

\mbox{}

\begin{pf*}{Proof of Proposition~\ref{Lemma2}}
Clearly $\be(v) \J\{v \in\CL(T_n)\}$ is $\CG_n$-measurable, and hence
$M_n(\gamma)$ is $\CG_n$-measurable.
We first prove that
%
\begin{equation}\label{martingale1}
\E[M_{n+1}{(\gamma)}|\CG_n]=M_n{(\gamma)}\bigl(1+\QQ(\gamma)/f_n\bigr).
\end{equation}
Given a sequence
$(T_n)_{n \geq1}$ of random $N$-ary recursive trees, one can define a
sequence $(V_n)_{n \geq1}$
of $\U$-valued random variables such that
\[
T_{n+1}= T_n \cup\{ (V_n,1),\ldots, (V_n,N)\}
\]
for every $n \geq0$, where $V_0=\varnothing$ and $V_n \in\CL(T_n)$.
The random variable $V_n$
corresponds to the random vertex chosen to generate $T_{n+1}$ from
$T_n$. Hence, by construction,
$ P(V_n=v| T_1,\ldots, T_n)= \J\{ v \in\CL(T_n)\} {1}/{f_n} $
for every $n\geq1$
and $P(V_n=v|\CG_n)=1/f_{n} \J\{v \in\CL(T_n)\} $. At this stage
one can
write
\begin{eqnarray*}
&&\E[M_{n+1}{(\gamma)}|\CG_n]\\
&&\qquad= M_n{(\gamma)}
+ \E\biggl[\sum_{v \in\CL(T_n)}
\be(v)^\gamma\bigl(A_1{(v)}^\gamma+\cdots+A_N{(v)}^\gamma-1
\bigr)\J\{
V_n=v\} \Big|\CG_n\biggr]\\
&&\qquad= M_n{(\gamma)} + \QQ(\gamma) \sum_{v \in\CL(T_n)} \be
(v)^\gamma\E[
\J\{V_n=v\}|\CG_n ]
=M_n{(\gamma)}\bigl(1+\QQ(\gamma)/f_n\bigr).
\end{eqnarray*}
Taking the expectation of both sides gives
$
\E[M_{n+1}{(\gamma)}]= \E[M_{n}{(\gamma)}](1+\QQ(\gamma)/f_n).
$
Since $\E[M_1{(\gamma)}]=\QQ(\gamma)+1$ and $f_0=1$, it follows
immediately that
%
\begin{equation}\label{40bis}
\E[M_{n}{(\gamma)}]=\prod_{i=0}^{n-1}\bigl(1+\QQ(\gamma
)/f_i\bigr)=m_n(\gamma).
\end{equation}
See (\ref{defmngamma}).
Since $f_i=(N-1)i+1$, by simple algebra one gets that
\[
m_n(\gamma)=\frac{\Gamma(({\CS(\gamma)+1})/({N-1})+n
)\Gamma
({1}/({N-1})) }{\Gamma(({\CS(\gamma)+1})/({N-1})
)\Gamma
({1}/({N-1})+n) }=\frac{(({\CS(\gamma
)+1})/({N-1})
)_n}{({1}/({N-1}))_n}.
\]
At this stage, recall that, given two positive real numbers $x$ and $y$,
\[
\frac{\Gamma(x+n)}{\Gamma(y+n)}=n^{x-y}\biggl(1+O\biggl(\frac
{1}{n}\biggr)\biggr)
\]
as $n \to+\infty$, which proves (\ref{asintm}).
Finally, (\ref{martingale1}) and (\ref{40bis}) yield that $\tilde
M_n{(\gamma)}$ is a $(\CG_n)_n$-martingale
since $M_n{(\gamma)} \geq0$ and $\E[M_n{(\gamma)}]<+\infty$ for every
$n\geq1$.
The last part of the theorem follows by classical martingale theory.
\end{pf*}
\begin{pf*}{Proof of Proposition~\ref{Lemma4}}
Observe that
\[
\bigl[\beta_{(n)}^{(\gamma)}\bigr]^\delta\leq\sum_{v \in\CL(T_n)} \frac
{\be
(v)^\delta}{ m_n(\gamma)^{{\delta}/{\gamma}}};
\]
hence for every $\eps>0$,
by Markov's inequality and (\ref{asintm}),
one gets
\begin{eqnarray*}
P\bigl\{ \b_{(n)}^{(\gamma)} >\eps\bigr\} & \leq &
P \biggl\{ \sum_{v \in\CL(T_n)} \frac{\be(v)^\delta}{m_n(\gamma
)^{
{\delta}/{\gamma}}} \geq\eps^\delta\biggr\}
\leq\frac{1}{\eps^\delta m_n(\gamma)^{\delta/\gamma}} \E
[M_n{(\delta
)} ]\\
&=& \frac{m_n(\delta) }{\eps^\delta m_n(\gamma)^{\delta/\gamma} }
\leq
\frac{ C_{\delta,\gamma}}{\eps^\delta}
n^{({\delta}/({N-1}))(\QQ(\delta)/\delta-\QQ(\gamma)/\gamma)}\\
&=&
\frac{ C_{\delta,\gamma}}{\eps^\delta}
n^{({\delta}/({N-1}))({\mu}(\delta)-{\mu}(\gamma))}.
\end{eqnarray*}
This proves the first statement.
When $\delta< \gamma$, one has $\delta/\gamma<1$, and hence, using
Minkowski inequality and (\ref{asintm}),
one gets
\begin{eqnarray*}
\E[\tilde M_n(\gamma)^{\delta/\gamma} ] & \leq & \frac{1}{
m_n(\gamma
)^{\delta/\gamma} }
\E\biggl[ \sum_{v \in\CL(T_n)} \be(v)^\delta\biggr] \\
&\leq& \frac{ m_n(\delta)}{ m_n(\gamma)^{\delta/\gamma} }
\leq C_{\delta,\gamma}n^{({\delta}/({N-1}))({\mu}(\delta)-{\mu
}(\gamma))},
\end{eqnarray*}
which proves the second statement. Assume now that $\delta>\gamma$.
In order to prove the last part of the statement
let us show that $\tilde M_n(\gamma)$ is uniformly integrable. To this
end, observe that
\[
\tilde M_{n+1}(\gamma)=\sum_{v \in\CL(T_n)}
\frac{\be(v)^\gamma}{m_n(\gamma)} \frac{[
1+(\sum_{k=1}^N A_k{(v)}^\gamma-1 )\J\{V_n=v\}
]}{1+
{\QQ(\gamma)}/{f_n}} ,
\]
and hence
\begin{eqnarray*}
\tilde M_{n+1}(\gamma)- \tilde M_{n}(\gamma)
&=&
- \frac{\tilde M_{n}(\gamma) \QQ(\gamma)}{f_n(1+\QQ(\gamma)/f_n)
}\\
&&{}+\frac{1}{m_{n+1}(\gamma)}\sum_{v \in\CL(T_n)} \be(v)^\gamma
\Biggl(\sum_{k=1}^N A_k{(v)}^\gamma-1 \Biggr)\J\{V_n=v\} .
\end{eqnarray*}
At this stage write
\begin{eqnarray*}
&&|\tilde M_{n+1}(\gamma) - \tilde M_{n}(\gamma)|^{{\delta
}/{\gamma
}}\\
&&\qquad\leq
2^{{\delta}/{\gamma}-1} \frac{|\tilde M_{n}(\gamma)|^{
{\delta
}/{\gamma}}|\QQ(\gamma)|^{{\delta}/{\gamma}} }{f_n^{
{\delta
}/{\gamma}} |1+\QQ(\gamma)/f_n|^{{\delta}/{\gamma}}}
\\
&&\qquad\quad{}
+\frac{2^{{\delta}/{\gamma}-1}}{ m_{n+1}(\gamma)^{
{\delta
}/{\gamma}}}
\sum_{v \in\CL(T_n)} \be(v)^{\delta}
\Biggl|\sum_{k=1}^N A_k{(v)}^\gamma-1 \Biggr|^{{\delta}/{\gamma
}} \J\{
V_n=v\}
\\
&&\qquad\leq
\frac{2^{{\delta}/{\gamma}-1}|\QQ(\gamma)|^{
{\delta
}/{\gamma}}}{|1
+\QQ(\gamma)/f_n|^{{\delta}/{\gamma}}}
\frac{1}{f_n^{{\delta}/{\gamma}}} f_n^{{{\delta}/{\gamma}}-1}
\sum_{v \in\CL(T_n)} \frac{\be(v)^{\delta}}{m_{n}(\gamma)^{
{\delta
}/{\gamma}}}
\\
&&\qquad\quad{}
+\frac{2^{{{\delta}/{\gamma}}-1}}{m_{n+1}(\gamma)^{
{\delta
}/{\gamma}}}
\sum_{v \in\CL(T_n)} \be(v)^{\delta}
\Biggl|\sum_{k=1}^N A_k{(v)}^\gamma-1 \Biggr|^{{\delta}/{\gamma
}}\J\{
V_n=v\} .
\end{eqnarray*}
Taking the expectation one gets
\begin{eqnarray*}
&&
\E|\tilde M_{n+1}(\gamma) - \tilde M_{n}(\gamma)|^{{\delta
}/{\gamma}}
\\
&&\qquad\leq
\frac{2^{{\delta}/{\gamma}-1}|\QQ(\gamma)|^{{\delta
}/{\gamma
}}}{|1+\QQ(\gamma)/f_n|^{{\delta}/{\gamma}}}
\frac{1}{f_n} \E\biggl[\sum_{v \in\CL(T_n)} \frac{\be(v)^{\delta
}}{m_{n}(\gamma)^{{\delta}/{\gamma}}} \biggr]
\\
&&\qquad\quad{}
+\frac{2^{{\delta}/{\gamma}-1}}{ m_{n+1}(\gamma)^{
{\delta
}/{\gamma}}}
\frac{1}{f_n} \E\biggl[ \sum_{v \in\CL(T_n)}\be(v)^{\delta} \biggr]
\E\Biggl|\sum_{k=1}^N A_k^\gamma-1 \Biggr|^{\delta/\gamma}
\\
&&\qquad\leq
C_1 \frac{1}{f_n}
\biggl[ \frac{ m_n( \delta)}{ m_n(\gamma)^{\delta/\gamma}}+ \frac{
m_n(\delta)}{ m_{n+1}(\gamma)^{\delta/\gamma}} \biggr]
\\
&&\qquad\quad   \mbox{by (\ref{asintm})}
\\
&&\qquad\leq C_2 \frac{1}{f_n} n^{{\delta({\mu}(\delta)-{\mu
}(\gamma))}/({N-1})}
\leq C_3 n^{{\delta({\mu}(\delta)-{\mu}(\gamma))}/({N-1}) -1}.
\end{eqnarray*}
Since ${\mu}(\delta)<{\mu}(\gamma)$, it follows that
%
\begin{equation}\label{convserie}
\sum_{i \geq1} \E[ |\tilde M_{i+1}(\gamma)-\tilde M_{i}(\gamma
)|^{
{\delta}/{\gamma}}] <+\infty.
\end{equation}
By the convexity of $\QQ(s)$, it is easy to see that
${\mu}(s)< {\mu}(\gamma)$ if $\gamma< s < \delta$.
Hence, without loss of generality, one can suppose that
$\gamma<\delta\leq2\gamma$.
Since $(\tilde M_n)_{n\geq1}$ is a martingale (cf. Proposition \ref
{Lemma2}) and $1<\delta/\gamma\leq2$, the
Topchii--Vatutin inequality (see, e.g., \cite{alsmeyerrosler1}) gives
\[
\E|\tilde M_{n}(\gamma)|^{{\delta}/{\gamma}} \leq\E|\tilde
M_1(\gamma)|^{{\delta}/{\gamma}}+ 2 \sum_{i=2}^{n}\E|\tilde
M_{i}(\gamma)-\tilde M_{i-1}(\gamma)|^{{\delta}/{\gamma}}.
\]
Combining this last
inequality with (\ref{convserie}), one obtains
\[
\sup_{n \geq1} \E|\tilde M_{n}(\gamma)|^{\delta/\gamma} <+\infty.
\]
Hence $(\tilde M_n(\gamma))_n$ is uniformly integrable and
then converges in $L^1$ to $\tilde M_\infty(\gamma)$ with $\E[\tilde
M_\infty(\gamma)]=\lim_n
\E[\tilde M_n(\gamma)]=1$.
\end{pf*}
\begin{pf*}{Proof of Proposition~\ref{propeqlimiteMinfty}}
Let $\psi_n(\xi)=\E[e^{i\xi\tilde M_n(\gamma)}]$.
Arguing as in the proof of Proposition~\ref{Propprobint}
and using the same notation, we get
\[
\tilde M_n(\gamma)=\sum_{j=1}^N A_j(\varnothing)^\gamma\frac
{m_{|T^{(j)}_n|}(\gamma)}{m_n(\gamma)}
\sum_{v \in\CL(T_{n}^{(j)})} \frac{\prod_{i=0}^{|v|-1} (
A_{v_{i+1}}^{(j)}(v|i))^\gamma}{ m_{|T^{(j)}_n|}(\gamma)}
\]
and then
\[
\psi_n(\xi)=\E\Biggl[ \prod_{j=1}^N
\psi_{|T^{(j)}_n|} \bigl( \xi A_j^\gamma\Delta_n^{(j)} \bigr)
\Biggr],
\]
where
\[
\Delta_n^{(j)}=\frac{m_{|T^{(j)}_n|}(\gamma)}{m_n(\gamma)}=
\biggl(\frac{|T^{(j)}_n|}{n} \biggr)^{{\QQ(\gamma)}/({N-1})}
\biggl(\frac{1+O({1}/{|T^{(j)}_n|})}{1+O(
{1}/{n}
)}\biggr).
\]
Now note that, for a suitable constant $C$,
$\Delta_n^{(j)}\leq C$ for every $j$ almost surely
and, by Proposition~\ref{comblemma}, $(\Delta_n^{(1)},\ldots,\Delta_n^{(N)})$
converges almost surely to $(U_1^{{\QQ(\gamma)}/{(N-1)}},\ldots
,U_N^{{\QQ(\gamma)}/{(N-1)}})$
where $(U_1,\ldots,U_N)$ has Dirichlet distribution
of parameters $(1/(N-1),\ldots,1/(N-1))$.
At this stage write
\[
\psi_n(\xi)=\E\Biggl[ \prod_{j=1}^N
\psi_{\infty,\gamma} \bigl( \xi A_j^\gamma\Delta_n^{(j)} \bigr)
\Biggr]+R_n
\]
with
\[
R_n=\E\Biggl[ \prod_{j=1}^N
\psi_{|T_n^{(j)}|} \bigl( \xi A_j^\gamma\Delta_n^{(j)} \bigr)
- \prod_{j=1}^N
\psi_{\infty,\gamma} \bigl( \xi A_j^\gamma\Delta_n^{(j)} \bigr)
\Biggr].
\]
By dominated convergence one gets
\[
\lim_{n \to+\infty} \E\Biggl[ \prod_{j=1}^N
\psi_{\infty,\gamma} \bigl( \xi A_j^\gamma\Delta_n^{(j)} \bigr)
\Biggr]=\E\Biggl[ \prod_{j=1}^N
\psi_{\infty,\gamma} \bigl( \xi A_j^\gamma U_j^{{\QQ(\gamma
)}/{(N-1)}} \bigr)
\Biggr].
\]
It remains to show that $R_n$ converges to zero. Recall that, given $2N$
complex numbers $a_1,\ldots,a_N,b_1,\ldots,b_N$
with $|a_i|, |b_i| \leq1$, then
\[
\Biggl|\prod_{i=1}^N a_i - \prod_{i=1}^N b_i\Biggr| \leq\sum_{i=1}^N |a_i-b_i|.
\]
Hence
\begin{eqnarray*}
|R_n| &\leq& \sum_{j=1}^N \E\bigl|\psi_{|T_n^{(j)}|} \bigl( \xi
A_j^\gamma\Delta_n^{(j)} \bigr)
- \psi_{\infty,\gamma} \bigl( \xi A_j^\gamma\Delta_n^{(j)}
\bigr)\bigr|
\\
&\leq& \sum_{j=1}^N \E\Bigl[ \sup_{x\dvtx |x| \leq|\xi A_j^\gamma C| }
\bigl|\psi_{|T_n^{(j)}|}(x)- \psi_{\infty,\gamma}(x)\bigr|\Bigr].
\end{eqnarray*}
Since point-wise convergence of characteristic functions yields
the same convergence on every compact set, and $\psi_n$
converges to $\psi_{\infty,\gamma}$, one has that
$\sup_{x\dvtx |x| \leq C }
|\psi_{n}(x)- \psi_{\infty,\gamma}(x)|$ converges to zero when $n$
goes to $+\infty$ for every $C<+\infty$.
By Proposition~\ref{comblemma} $|T_n^{(j)}|$ converges almost surely to
$+\infty$, hence
dominated convergence yields that
$ \sup_{x\dvtx |x| \leq|\xi A_j^\gamma C| }
|\psi_{|T_n^{(j)}|}(x)- \psi_{\infty,\gamma}(x)| $
converges almost surely to zero as $n$ goes to $+\infty$ and then, by
dominated convergence,
$R_n$ converges to zero.
\end{pf*}

In order to prove Proposition~\ref{proplimit1}, let us consider the
following central limit result.
Assume that
$(a_{jn})_{jn}$ is an array of positive weights, and let $r_n$ be a
diverging sequence of integer numbers.
Given any array of identically distributed and row-wise independent
random variables $(X_{jn})_{n \geq1, j=1,\ldots,r_n}$ with probability
distribution function $F_0$,
set
\[
\tilde S_n:=\sum_{j=1}^{r_n} a_{jn} X_{jn}.
\]
Moreover assume that, for some $\gamma$ in $(0,2]$,
%
\begin{equation}\label{H0}
\lim_{n \to+\infty} \sum_{j=1}^{r_n} a_{jn}^\gamma=a_\infty
\quad \mbox{and}\quad  \lim_{n \to+\infty} \max_{j=1,\ldots,r_n} a_{jn}
= 0.
\end{equation}

It is not hard to prove the following.
%
\begin{lemma}\label{Lemma5}
Let (\ref{H0}) and~\ref{hypoHgama} be in force for $\gamma$
in $(0,2]$.
Then
%
\begin{equation}
\label{eqchar02bis}
\lim_{n \to+\infty} \E[ e^{ i \xi\tilde S_{n}}]
=\hat g_\gamma(\xi a^{{1}/{\gamma}}_{\infty})
\end{equation}
for every $\xi\in\RE$, $\hat g_\gamma$ being defined in (\ref{chaSta}).
\end{lemma}
\begin{pf} The proof can be obtained, following the same line of
the proofs of Lemmas 4, 5, 6 in \cite{BLM2008}, as a consequence of the
central limit
theorem for triangular array; see, for example, \cite{fristedgray}.
\end{pf}
\begin{pf*}{Proof of Proposition~\ref{proplimit1}}
Denote by $\CB$ the $\s$-algebra generated by $\{T_n,\beta_{j,n}\dvtx n
\geq1, j=1,\ldots, f_n \}$.
The proof is
essentially an application of Lem\-ma~\ref{Lemma5}
to the conditional law of $ \tilde W_n $ given $\CB$.
By Propositions~\ref{Lemma2} and~\ref{Lemma4}, every divergent sequence
$(n')$ of integer numbers
contains a divergent subsequence $(n'') \subset(n')$
for which $\tilde M_{n''}(\gamma)$ converges almost surely to $\tilde
M_\infty(\gamma)$,
and $\beta_{(n'')}^{(\gamma)}$ converges almost surely to zero.
Hence by Lemma~\ref{Lemma5} we obtain
$
\lim_{n'' \to+\infty}\E[e^{i\xi\tilde W_{n''}}|\CB]=\hat g_\gamma
(\xi
\tilde M_{\infty}^{{1}/{\gamma}}(\gamma))
$
almost surely. Dominated convergence theorem yields
$
\lim_{n'' \to+\infty}\E[e^{i\xi\tilde W_{n''}}]=\E[\hat g_\gamma
(\xi
\tilde M_{\infty}^{{1}/{\gamma}}(\gamma))]
$.
Since the limiting function is independent of the arbitrarily chosen
sequence $(n')$,
a classical argument shows that the last limit is true with $n\to
+\infty$ in place of
$n''\to+\infty$.
\end{pf*}
\begin{pf*}{Proof of Proposition~\ref{Propjointconv}}
Let us first prove that
when $t$ goes to $+\infty$,
$\nu_t e^{-t(N-1)}$ converges in distribution
to a random variable $\GA$ with Gamma distribution of parameters $(1/(N-1),1)$.
Since $\nu_t$ is a negative binomial random variable of parameters
$(1/(N-1),\exp\{-(N-1)t\})$,
for every integer~$k$
\[
P\{ \nu_t \leq k\}=
\frac{\Gamma(k+1+{1}/({N-1}))}{\Gamma(k+1)\Gamma(1/(N-1))}
\int_0^{e^{-(N-1)t}} x^{{1}/({N-1})-1}(1-x)^{k}\,dx;
\]
see formula (5.31) in \cite{JohnsonKotzKempUnidisc}.
Hence, for every $y>0$, after setting $k_t=\lfloor y
e^{(N-1)t}\rfloor$
(where $\lfloor x \rfloor$ is the integer part of $x$),
one can write
\begin{eqnarray*}
&&
P\bigl\{ \nu_t e^{-(N-1)t} \leq y\bigr\}\\
&&\qquad= P\{ \nu_t \leq k_t\} \\
&&\qquad= \frac{\Gamma(k_t+1+{1}/({N-1}))}{\Gamma(k_t+1)\Gamma(1/(N-1))}
\int_0^{e^{-(N-1)t}} x^{{1}/({N-1})-1}(1-x)^{k_t}\,dx \\
&&\qquad= \frac{\Gamma(k_t+1+{1}/({N-1}))}{\Gamma(k_t+1)\Gamma(1/(N-1))}
\frac{1}{[ye^{(N-1)t}]^{{1}/({N-1})}}\\
&&\qquad\quad{}\times\int_0^y u^{{1}/({N-1})-1}
\biggl(1-\frac{u}{ye^{(N-1)t}}\biggr)^{k_t}\,du\\
&&\qquad= \frac{1+O(1/k_t)}{\Gamma(1/(N-1))}
\biggl[\frac{k_t}{ye^{(N-1)t}}\biggr]^{{1}/({N-1})}\\
&&\qquad\quad{}\times\int_0^y u^{{1}/({N-1})-1} \biggl(1-\frac{u}{ye^{(N-1)t}}
\biggr)^{ye^{(N-1)t} {k_t}/({ye^{(N-1)t}})}\,du.
\end{eqnarray*}
Since ${k_t}/{ye^{(N-1)t}} \to1$, by dominated convergence one gets
\[
\lim_{t \to+\infty}
P\bigl\{ \nu_t e^{-(N-1)t} \leq y\bigr\}=\frac{1}{\Gamma({1}/({N-1}))} \int_0^y
u^{{1}/({N-1})-1}e^{-u}\,du.
\]
At this stage, since $\nu_t$ converges in probability to $+\infty$,
(\ref{asintm}),
Slustky's theorem and the continuous mapping theorem yield that
%
\begin{equation}\label{rr2}
\lim_{t \to+\infty}\E\bigl[e^{i\xi N_t(\gamma)}\bigr]=\E\bigl[e^{i\xi c_\gamma
\GA
^{{{\mu}(\gamma)}/{(N-1)}} } \bigr]
\end{equation}
for every $\xi$ in $\RE$. Setting $u_{n}(\xi):=\E[e^{i\xi\tilde
W_n}]$ by
Proposition~\ref{proplimit1} we know that
%
\begin{equation}\label{rr3}
\lim_{n \to+\infty}u_{n}(\xi) = \E[\hat g_\gamma(\xi\tilde
M_\infty
(\gamma)^{{1}/{\gamma}})]=: u_{\infty}(\xi)
\end{equation}
for every $\xi$ in $\RE$.
For every diverging sequence $(t_n)_n$
write
\[
\E\bigl[e^{i\xi_1 N_{t_{n}}(\gamma) +i\xi_2 \tilde W_{\nu_{t_{n}}}}\bigr]
=u_{\infty}(\xi_2) \E\bigl[e^{i\xi_1 N_{t_{n}}(\gamma)}\bigr]+R_{n},
\]
where
\[
R_{n}=\E\bigl[e^{i\xi_1 N_{t_{n}}(\gamma)}\bigl(e^{i\xi_2 \tilde W_{\nu
_{t_{n}}}}-u_{\infty}(\xi_2)\bigr)\bigr].
\]
It is easy to show that
\[
\lim_{n\to+\infty} |R_n| \leq\lim_{n\to+\infty} \E|u_{\nu
_{t_n}}(\xi
_2)- u_{\infty}(\xi_2)|=0
\]
by dominated convergence, since
$\nu_{t_{n}}$ converges in probability to $+\infty$, and (\ref{rr3}) holds.
The result now follows from (\ref{rr2}).
The second part of the statement follows immediately by the continuous
mapping theorem.
\end{pf*}

\subsection{\texorpdfstring{Proofs of Section \protect\ref{SecMain}}{Proofs of Section 2}}\label{proofsMain}

Under the hypotheses of Proposition~\ref{proplimit1},
(\ref{limitw}) yields that, when $\gamma\not=1$ or when $\gamma=1$ and
condition (b) of~\ref{hypoH1} holds,
$e^{-\mu(\gamma)t}V_t$ converges in distribution
to a scale mixture of stable laws. The scale mixing measure is the law
of $c_\gamma\GA^{{\mu(\gamma)}/({N-1})}
\tilde M_\infty(\gamma)^{{1}/{\gamma}}$, with $\GA$ and
$\tilde
M_\infty(\gamma)$ stochastically
independent. While if $\gamma=1$ and condition (a) of~\ref{hypoH1} holds,
then $e^{-\mu(\gamma)t}V_t$ converges in distribution to $m_0 c_1 \GA
^{{\mu(1)}/({N-1})}
\tilde M_\infty(1)$.
Again the mixing measure is the law of $c_1 \GA^{{\mu(1)}/({N-1})}
\tilde M_\infty(1)$.
At this stage, in order to complete the proof of the main theorems of
Section~\ref{SecMain1}, it remains to
study in more detail the distribution of $c_\gamma\GA^{{\mu
(\gamma
)}/({N-1})}\tilde M_\infty(\gamma)^{{1}/{\gamma}}$.

The more interesting case is
${\mu}(\delta)<{\mu}(\gamma)$ for $\delta>\gamma$.
Proposition~\ref{propeqlimiteMinfty} shows that the law of
$\tilde M_\infty(\gamma)$ satisfies the fixed\vspace*{1pt} point equation
for distributions~(\ref{integraleqpsi2}).
We will show that the law of $c^\gamma_\gamma\GA^{{\QQ(\gamma)}/({N-1})}
\tilde M_\infty(\gamma)$ and
the limit law of $e^{-\mu(\gamma)t}V_t$ satisfy a fixed point
equation too. In view of known results on this kind of equations,
we will be able to complete the proof of Theorem~\ref{teo2} of Section
\ref{SecMain1}.

In what follows denote by $\operatorname{Beta}(a,b)$ [$\operatorname{Gamma}(a,b)$, resp.],
$a>0$ and $b>0$, the beta distribution of parameters
$a$ and $b$ (the gamma distribution of shape parameter $a$ and scale
parameter $b$, resp.).
We will need the following result.

\begin{lemma}\label{Lemma8}
Let $\GA_1,\ldots,\GA_N,V$ be independent random variables
where~$V$ has $\operatorname{Beta}(1/(N-1),1)$ distribution, and
$\GA_i$ has $\operatorname{Gamma}(1/(N-1),1)$ distribution for every $i$. Let $\GA$
and $U=(U_1,\ldots,U_N)$ be stochastically independent, $\GA$~with
$\operatorname{Gamma}(1/(N-1),1)$
distribution and $U$ with
Dirichlet distribution of parameters $(1/(N-1),\ldots,1/(N-1))$. Then
\[
(\GA U_1,\GA U_2,\ldots,\GA U_N) \stackrel{\CL} = (V\GA_1,V\GA
_2,\ldots
,V\GA_N).
\]
\end{lemma}
\begin{pf}
Set $S=\sum_{i=1}^N \GA_i$. $S$ is a~$\operatorname{Gamma}(N/(N-1),1)$ random
variable, and $S$ and $V$
are stochastically independent.
It is easily seen that $SV$ is a~$\operatorname{Gamma}(1/(N-1),1)$ random variable.
Now define $\tilde U:=(\GA_1/S,\ldots,\GA_N/S)$ and $\tilde\GA
:=SV$. It
is well known that
$\tilde U$ has a Dirichlet distribution of parameters $(1/(N-1),\ldots
,1/(N-1)$, and
$S$ and $\tilde U$ are independent; see, for example, Section 10.4 in
\cite{fristedgray}.
Hence, it turns out that $\tilde U$ and $\tilde\GA$ are stochastically
independent. Clearly
\[
(\tilde\GA\tilde U_1,\ldots,\tilde\GA\tilde U_N)=(V\GA_1,\ldots
,V\GA_N),
\]
which proves the claim.
\end{pf}
%
\begin{proposition}\label{Lemma7bis} Let the assumptions of
Proposition~\ref{proplimit1}
be in force, and let $v_{\infty,\gamma}$ be the
characteristic function of $c^\gamma_\gamma\GA^{{\QQ(\gamma)}/({N-1})}
\tilde M_\infty(\gamma)$. Then~$v_{\infty,\gamma}$ satisfies
the integral equation (\ref{eqv}), that is,
%
\begin{equation}\label{inteqmix}
v_{\infty,\gamma}(\xi)=\int_0^1
\E\Biggl[\prod_{i=1}^N v_{\infty,\gamma}\bigl(\xi A_i^\gamma\tau^{\QQ
(\gamma
)}\bigr) \Biggr] \,d\tau.
\end{equation}
Moreover, if $w_{\infty,\gamma}$
denotes the characteristic function
of the limit in distribution of $e^{-\mu(\gamma)t}V_{t}$,
then $w_{\infty,\gamma}$
satisfies equation (\ref{stazselfintegral}) for $\mu^*=\mu(\gamma)$,
that is,
%
\begin{equation}\label{integraleqw}
w_{\infty,\gamma}(\xi)= \int_{0}^1
\E\Biggl[\prod_{i=1}^N w_{\infty,\gamma}\bigl(A_i \tau^{\mu(\gamma)}
\xi\bigr)
\Biggr]\,d\tau
=\int_{0}^1 \hat{Q}[w_{\infty,\gamma}]\bigl(\tau^{\mu(\gamma)}
\xi\bigr)\,d\tau.\hspace*{-35pt}
\end{equation}
\end{proposition}
\begin{pf}
Recall that $\psi_{\infty,\gamma}$ denotes the characteristic
function of
$\tilde M_\infty(\gamma)$. Hence, from the independence of
$\GA$ and $\tilde M_\infty(\gamma)$,
we have
\[
v_{\infty,\gamma}(\xi)=\E\bigl[\psi_{\infty,\gamma}\bigl(\xi c^\gamma
_\gamma\GA
^{{\QQ(\gamma)}/({N-1})} \bigr)\bigr].
\]
Since $\psi_{\infty,\gamma}$ satisfies equation (\ref{integraleqpsi}),
we can write
\[
v_{\infty,\gamma}(\xi)=\E\Biggl[\prod_{i=1}^N \psi_{\infty,\gamma
}\bigl(A_i^{\gamma}(U_i\GA)^{{\QQ(\gamma)}/({N-1})}
c_\gamma^\gamma\xi\bigr)\Biggr],
\]
where $U=(U_1,\ldots,U_N)$, $(A_1,\ldots,A_N)$ and $\GA$ are independent,
$U$ has Dirichlet distribution of parameters $(1/(N-1),\ldots,1/(N-1))$
and $\GA$ has\break $\operatorname{Gamma}(1/(N-1),1)$ distribution.
By Lemma~\ref{Lemma8} if
$(\GA_1,\ldots,\GA_N,V)$ are independent random variables,
$\GA_i$ with $\operatorname{Gamma}(1/(N-1),1)$ distribution, $V$ with
$\operatorname{Beta}(1/(N-1),1)$ distribution
and $(\GA_1,\ldots,\GA_N,V)$ and $(A_1,\ldots,A_N)$ independent,
then we can write
\begin{eqnarray*}
v_{\infty,\gamma}(\xi)&=&
\E\Biggl[\prod_{i=1}^N \psi_{\infty,\gamma}\bigl(A_i^{\gamma}(V\GA
_i)^{
{\QQ(\gamma)}/({N-1})} c_\gamma^\gamma\xi\bigr)\Biggr]\\
&=&\E\Biggl[\prod_{i=1}^N v_{\infty,\gamma}\bigl(A_i^{\gamma}V^{{\QQ
(\gamma
)}/({N-1})} \xi\bigr)\Biggr].
\end{eqnarray*}
Then (\ref{inteqmix}) follows since $V^{1/(N-1)}$ has uniform
distribution on $(0,1)$.

As for the second part, let $\hat g_\gamma$ be defined in (\ref{chaSta}).
From Proposition~\ref{Propjointconv}
we know that
\[
w_{\infty,\gamma}(\xi)=\E\bigl[ \hat g_\gamma\bigl(c_\gamma\GA^{
{\mu
(\gamma)}/({N-1})}
\tilde M_\infty(\gamma)^{{1}/{\gamma}} \bigr) \bigr].
\]
Define
$ \YY=c^\gamma_\gamma\GA^{{\QQ(\gamma)}/({N-1})}
\tilde M_\infty(\gamma) $; then (\ref{inteqmix}) is equivalent to
%
\begin{equation}\label{inteqmix2}
\YY\stackrel{\CL}= \Theta^{\QQ(\gamma)} \sum_{i=1}^N A_i^\gamma
\YY_i,
\end{equation}
where $(\YY,\YY_1,\ldots,\YY_N)$ are i.i.d., $\Theta$ has uniform
distribution on $(0,1)$
and $(\YY,\YY_1,\ldots,\YY_N)$, $\Theta$ and $(A_1,\ldots,A_N)$
are independent. Then from (\ref{inteqmix2}) and from the analytic form
of $\hat g_\gamma$ we get
\begin{eqnarray*}
w_{\infty,\gamma}(\xi) & = &
\E[ \hat g_\gamma(\xi\YY^{{1}/{\gamma}}) ]
=\E\Biggl[ \prod_{i=1}^N \hat g_\gamma\bigl(\xi\Theta^{{\mu}(\gamma)}A_i
\YY_i^{{1}/{\gamma}}\bigr) \Biggr] \\
& = & \E\Biggl[ \prod_{i=1}^N w_{\infty,\gamma}\bigl( \xi\Theta^{{\mu
}(\gamma)}A_i\bigr)
\Biggr],
\end{eqnarray*}
and this completes the proof.
\end{pf}

In order to prove Proposition~\ref{Pfixpoint} we need to recall few
important results on fixed point equations
for distributions.
Assume that
$B=(B_1,\ldots,B_N)$ is a vector of nonnegative random variables.
Consider the following fixed point equation:
%
\begin{equation}\label{fixedEq2}
\nu=T(\nu),
\end{equation}
where, given any probability distribution $\nu$, $T(\nu)$ is the law
of $\sum_{i=1}^N B_i Y_i$, where $B$ and $(Y_1,\ldots,Y_N)$
are stochastically independent, and $(Y_1,\ldots,Y_N)$ are stochastically
independent and identically distributed random variables
with distribution~$\nu$. Clearly, (\ref{fixedEq2}) is equivalent to
equation
%
\begin{equation}\label{eqthesCF}
\phi(\xi)=\E\Biggl[\prod_{i=1}^N \phi(B_i\xi) \Biggr]
\end{equation}
for the corresponding Fourier--Stieltjes transforms.
Equations (\ref{eqv}), (\ref{integraleqpsi}), (\ref{inteqmix}) and
(\ref{integraleqw})
have this form for a suitable $B$.
In order to describe the fixed points of (\ref{fixedEq2}),
we introduce the convex function $a\dvtx[0,\infty)\to[0,\infty]$~by
%
\begin{equation}
\label{eqthesBis}
a(s):=\E\Biggl[\sum_{j=1}^N B_j^s \Biggr]
\end{equation}
with the convention that $0^0=0$.
%
\begin{proposition}[(\cite{DurrettLiggett1983,Liu1998,Liu2000,alsmeyerbiggins})]
\label{Prop0-1}
Assume that condition (\ref{H0-1}) holds true with
$B_i$ in place of $A_i$, that $P\{B_i=0 \mbox{ or }
1 \ \forall i=1,\ldots,N\}<1$ and that $a(1)=1$.
\begin{longlist}[(ii)]
\item[(i)] If $\sum_{j=1}^N B_j =1$ almost surely, 
then $a (s) \geq1$ for every $s<1$ and $a (s) \leq1$ for every $s>1$.
Moreover, the unique solution $\zeta$ of (\ref{fixedEq2})
with $\int_{\RE^+} \zeta(dv)=\int_{\RE^+} v \zeta(dv)=1$
is
the degenerate probability distribution $\zeta(\cdot)=\delta_1(\cdot)$;
\item[(ii)] If $P\{\sum_{j=1}^N B_j =1\}<1$
and if $a(s)<1$ for some $s>1$,
then (\ref{fixedEq2}) has a unique solution $\zeta$ with
$\int_{\RE^+}\zeta(dv)=\int_{\RE^+} v \zeta(dv)=1$.
Moreover $\zeta$ is nondegenerate and, for any $p>1$,
$\int_{\RE^+} v^{p} \zeta(dv)<+\infty$ if and only if $a(p)<1$.
\end{longlist}
\end{proposition}
\begin{pf*}{Proof of Proposition~\ref{Pfixpoint}}
Equation (\ref{eqv}) is of type (\ref{eqthesCF})
with $B_i=A_i^\gamma\Th^{\QQ(\gamma)}$, $\Th$ being
an uniform random variable on $[0,1]$ independent from $(A_1,\ldots,A_N)$.
Hence, in this case, the function $a$ defined in (\ref{eqthesBis}) becomes
\[
a(s):=a_\gamma(s)=\frac{\QQ(\gamma s)+1}{\QQ(\gamma) s+1}.
\]
Clearly $a_\gamma(1)=1$. Now, since $\delta>\gamma$
and $\QQ(\gamma)/\gamma={\mu}(\gamma)>{\mu}(\delta)\geq
-1/\delta$,
it is easy to see that the convex function $q \mapsto a_\gamma
(q/\gamma)$
is equal to $1$ in $q=\gamma$
and strictly smaller than $1$ in $q=\delta$. Since ${\mu}(q)-{\mu
}(\gamma)=0$
if and only if $a_\gamma(q/\gamma)=1$, it follows that\vadjust{\goodbreak}
equation ${\mu}(q)-{\mu}(\gamma)=0$ has at most one solution
$q_\gamma^*\not
=\gamma$.
This proves (ii). Noticing that $\delta/\gamma>1$ and
$a_\gamma(\delta/\gamma)<1$, by Proposition~\ref{Prop0-1},
(i) follows.
Since $\Th$ and $A$ are independent, $\Th^{\QQ(\gamma)}\sum_{i=1}^N
A_i^\gamma=1$
almost surely if and only if $\sum_{i=1}^N A_i^\gamma=1$ almost surely.
Hence, by (ii) of Proposition~\ref{Prop0-1}, $\zeta_{\infty,\gamma}$
is degenerate if and only if $\sum_{i=1}^N A_i^\gamma=1$ almost surely.
Finally, using that $a_\gamma$ is convex, $a_\gamma(1)=1$ and
$a_\gamma(\delta/\gamma)<1$, it follows that
$a_\gamma(p/\gamma)<1$ if and only if $p<q_\gamma^*$. Again by (ii)
of Proposition~\ref{Prop0-1}, whenever $\zeta_{\infty,\gamma}$
is nondegenerate,
$\int_{\RE^+} v^{{p}/{\gamma}} \zeta_{\infty,\gamma
}(dv)<+\infty$
if and only $p<q_\gamma^*$.
\end{pf*}
\begin{pf*}{Proof of Theorem~\ref{teo2}}
Proposition~\ref{Lemma4} yields that
$\E[\tilde M_\infty(\gamma)]=1$, since ${\mu}(\gamma)>{\mu
}(\delta)$ for
$\delta> \gamma$.
An easy computation shows that
$\E[c_\gamma^\gamma\GA^{{{\mu}(\gamma)\gamma}/({N-1})}]=1$
and then
$\E[c_\gamma^\gamma\GA^{{{\mu}(\gamma)\gamma}/({N-1})}\tilde
M_\infty
(\gamma)]=1$.
Hence, recalling that
$v_{\infty,\gamma}$ is the characteristic function of
$c_\gamma^\gamma\GA^{{{\mu}(\gamma)\gamma}/({N-1})} \tilde
M_\infty
(\gamma)$,
by (\ref{inteqmix}) of Proposition~\ref{Lemma7bis} and (i) of
Proposition~\ref{Pfixpoint},
the law of
$c_\gamma^\gamma\GA^{{{\mu}(\gamma)\gamma}/({N-1})}\tilde
M_\infty
(\gamma)$
is equal to $\zeta_{\infty,\gamma}$. At this stage (\ref{lim-teo2})
follows by (\ref{limitw}).
Moreover, from (\ref{integraleqw}) of Proposition~\ref{Lemma7bis},
$w_{\infty,\gamma}$ is a solution of (\ref{stazselfintegral})
for $\mu^*={\mu}(\gamma)$.
The proof of (i) is complete.

In order to prove (ii) let us observe that, from the properties of
$\gamma$-stable distributions, it follows that
$\int_{\RE^+} v^p \rho_{\infty,\gamma}(dv)<+\infty$
if and only if $p<\gamma$ and $\int_{\RE^+} v^{{p}/{\gamma}}
\zeta
_{\infty,\gamma}(dv)<+\infty$,
but for $p<\gamma$,
\[
\int_{\RE^+} v^{{p}/{\gamma}} \zeta_{\infty,\gamma}(dv) \leq
\biggl(\int_{\RE^+} v \zeta_{\infty,\gamma}(dv)\biggr)^{
{p}/{\gamma}}=1.
\]
It remains to show that $\rho_{\infty,\gamma}$ is a $\gamma$-stable
distribution
if and only if $\sum_{i=1}^N A_i^\gamma=1$ almost surely. This
follows from
(iii) of Proposition~\ref{Pfixpoint} and the fact that
%
\begin{eqnarray}\label{mixturestable}
&&
e^{ -|\xi|^\gamma k_1 (1-i \eta_1 \tan(\pi\gamma/2) \operatorname
{sign}\xi)} \\
&&\qquad  =\int_{\RE^+} e^{ -|\xi|^\gamma z k_0 (1-i \eta_0
\tan(\pi\gamma/2)\operatorname{sign}\xi) }
\zeta_{\infty,\gamma}(dz)
\nonumber
\end{eqnarray}
if and only if $k_1=k_0$, $\eta_1=\eta_0$ and $\zeta_{\infty,\gamma
}=\delta_1$.
Let us prove the last claim. Write (\ref{mixturestable}) for $\xi>0$ with
$\xi^\gamma=x$, and differentiate the resulting identity with respect
to $x$
to obtain
%
\begin{eqnarray}\label{mixturestable2}
&&
- k_1 \bigl(1-i \eta_1 \tan(\pi\gamma/2) \bigr) e^{ -x k_1 (1-i \eta_1 \tan
(\pi
\gamma/2) )} \nonumber\\[-8pt]\\[-8pt]
&&\qquad =\int_{\RE^+} k_0 z \bigl(1-i \eta_0
\tan(\pi\gamma/2)\bigr) e^{ -x k_0 (1-i \eta_0
\tan(\pi\gamma/2))}
\zeta_{\infty,\gamma}(dz).
\nonumber
\end{eqnarray}
Taking the limit
for $x \downarrow0$, recalling that $\int_{\RE^+} z \zeta_{\infty
,\gamma}(dz)=1$,
by dominated convergence
one gets
\[
k_1 \bigl(1-i \eta_1 \tan(\pi\gamma/2) \bigr)=k_0 \bigl(1-i \eta_0 \tan(\pi
\gamma/2) \bigr)
\]
and hence $k_1=k_0$ and
$\eta_0=\eta_1$. At this stage it suffices to
recall that a~scale mixture of stable laws is an identifiable family of
distributions.
See, for example, \cite{Teicher61}.

Analogously, (iii) and (iv) follow from (\ref{inteqmix}) of
Proposition~\ref{Lemma7bis}
and from~(iii) of Proposition~\ref{Pfixpoint}.
\end{pf*}
\begin{pf*}{Proof of Theorem~\ref{teo1}}
Recall that ${\mu}(\delta)<{\mu}(\gamma)$ for $\delta<\gamma$,
hence by Proposition
\ref{Lemma4} yields that $\tilde M_\infty(\gamma)=0$, and this
completes the proof.
\end{pf*}
\begin{pf*}{Proof of Theorem~\ref{PropW-2}}
We shall assume that $\W_\delta(X_0,V_\infty)<+\infty$, since otherwise
the claim is trivial.
Then, there exists an optimal pair $(X^*,Y^*)$ realizing the infimum
in the definition of the Wasserstein distance,
%
\begin{equation}
\label{eqdefdel}
\Delta:=\W_{\delta}^{\max(\delta,1)}(X_0,V_{\infty})=\W_\delta
^{\max
(\delta,1)}(X^*,Y^*)=\E|X^* -Y^* |^\delta.
\end{equation}
Let $(X_v^*,Y_v^*)_{v \in\U}$ be a sequence of independent and
identically distributed random variables
with the same law of $(X^*,Y^*)$,
which are further independent of $(\nu_t)_{t \geq0}$, $(T_n)_{n \geq
1}$, $(A(v))_{v \in\U}$.
By Proposition~\ref{Propprobint} it follows that $\sum_{j=1}^{f_{\nu
_t}} X_{j,\nu_t}^* \b_{j,\nu_t}$ has the same law of $V_t$,
where $X^*_{j,n}=X^*_{L_{j,n}}$ and $L_{j,n}$ is defined at the end of
Section~\ref{Srandtreeandprobrep}.
Moreover, since the characteristic function of~$V_\infty$ is a
solution of
(\ref{stazselfintegral}) with \mbox{$\mu^*=\mu(\gamma)$}, as already noted in
the \hyperref[intro]{Introduction},
the characteristic function of $e^{\mu(\gamma)t}V_\infty$ is a
solution of
(\ref{eqboltzivp}) with $\phi_0=w_{\infty,\gamma}$. Hence, applying
once again
Proposition~\ref{Propprobint}, we get that $e^{\mu(\gamma
)t}V_\infty$
has the same law of
$
\sum_{j=1}^{f_{\nu_t}} Y_{j,\nu_t}^* \b_{j,\nu_t}
$, where $Y^*_{j,n}=Y^*_{L_{j,n}}$.
For the sake of simplicity write $(X_j^*,Y^*_j)$ in place of
$(X_{j,n}^*,Y_{j,n}^*)$.
We can write
\begin{eqnarray*}
&&
\W_\delta^{\max(\delta,1)}
\bigl(e^{-\mu(\gamma)t}V_t,V_{\infty}\bigr)\\
&&\qquad = \W_\delta^{\max(\delta,1)} \bigl(e^{-\mu(\gamma)t}V_t,e^{-\mu
(\gamma
)t}e^{\mu(\gamma)t}V_{\infty}\bigr) \\
&&\qquad= e^{-\delta\mu(\gamma)t} \W_\delta^{\max(\delta,1)}
\bigl(V_t,e^{\mu
(\gamma)t}V_{\infty}\bigr) \\
&&\qquad \leq e^{-\delta\mu(\gamma)t} \sum_{n \geq0} \zeta(t,n)
\E\Biggl[ \Biggl|\sum_{j=1}^{f_n} X_{j}^* \b_{j,n} - \sum_{j=1}^{f_n}
Y_{j}^* \b_{j,n} \Biggr|^\delta\Biggr] \\
&&\qquad= e^{-\delta\mu(\gamma)t} \sum_{n \geq0} \zeta(t,n)\E\Biggl[ \E
\Biggl[ \Biggl| \sum_{j=1}^{f_n} (X_{j}^*-Y_{j}^*) \b_{j,n} \Biggr|^\delta
\bigg|\CG_n \Biggr] \Biggr],
\end{eqnarray*}
where $ \zeta(t,n)$ is the density of $\nu_t$ [see (\ref
{defnegativebinomial})] and $\CG_n=\s(A(v)_{v \in\CI
(T_n)},T_1,\ldots,T_n)$.
Now, if $0<\gamma<\delta\leq1$,
then Minkowski's inequality yields
%
\begin{eqnarray}\label{eqwascaso1}\quad
\E\Biggl[ \E\Biggl[ \Biggl| \sum_{j=1}^{f_n} (X_{j}^*-Y_j^*) \b_{j,n}
\Biggr|^\delta\bigg|\CG_n \Biggr] \Biggr]
&\leq& \E\Biggl[ \E\Biggl[ \sum_{j=1}^{f_n} \b_{j,n}^\delta
|X_j^*-Y_j^*|^\delta\bigg| \CG_n \Biggr] \Biggr]
\nonumber\\[-8pt]\\[-8pt]
& = &\E\Biggl[ \sum_{j=1}^{f_n} \b_{j,n}^\delta\Biggr] \Delta,
\nonumber
\end{eqnarray}
where $\Delta$ is defined in (\ref{eqdefdel}).
We now want to prove a similar inequality for $1\leq\gamma< \delta
\leq2$.
First of all we need to observe that, in addition to $\E
|X^*_j-Y_j^*|^\delta= \W_\delta^\delta(X_0,V_\infty)<+\infty$,
we have also that $\E(X^*_{j}-Y_{j}^*)=0$. If $\gamma\not=1$ the
claim follows since
by hypothesis $\E(X^*_{j})=\E(X_0)=0$ and $\E( Y_{j}^*)=\E(
V_\infty
)=0$, thanks to the fact that
$V_\infty$ is a mixture of centered stable random variables of
exponent $\gamma>1$.
When $\gamma=1$ and (a) of~\ref{hypoH1} holds, the proof of
the claim is similar.
When $\gamma=1$ and (b) of~\ref{hypoH1} holds, the proof
requires more care, since
$\E|X_0|=\E|V_\infty|=+\infty$.
Let $F_\infty(y)$ be the probability distribution function of~$V_\infty$, that is,
$F_\infty(x)=\int_{(-\infty,x]} \rho_{\infty,\gamma}(dy)$,
and recall that as an optimal pair one can choose
$(X^*,Y^*)=(F_0^{-1}(U),F_\infty^{-1}(U))$, $U$ being
a uniform random variable on $(0,1)$ and $F_0^{-1}$ ($F_\infty^{-1}$,
resp.) is the quantile function of
$F_0$ ($F_\infty$, resp.); see, for example, \cite{Rachev1991}.
Note that $\E|X^*_j-Y_j^*|^\delta<+\infty$, which yields that
$\E|X^*_j-Y_j^*|=\int_0^1|F_0^{-1}(u)-F_\infty^{-1}(u)|\,du <+\infty$.
Since $F_0$ and $F_\infty$ are symmetric distribution functions,
it is easy to see that $F_0^{-1}(U)-F_\infty^{-1}(U)$ is a symmetric
random variable, and hence
$\E(F_0^{-1}(U)-F_\infty^{-1}(U))=\E(X^*_{j}-Y_{j}^*)=0$.
Summarizing, if $1 \leq\gamma\leq2$, we have $\E(X^*_{j}-Y_{j}^*)=0$
and $\E|X^*_j-Y_j^*|^\delta<+\infty$. Hence
we can apply the Bahr--Esseen inequality (see \cite{BahrEsseen1965})
to obtain\looseness=-1
%
\begin{eqnarray}\label{eqwascaso2}
\E\Biggl[ \E\Biggl[ \Biggl| \sum_{j=1}^{f_n} (X_{j}^*-Y_j^*) \b_{j,n}
\Biggr|^\delta\bigg| \CG_n \Biggr] \Biggr]
&\leq& \E \Biggl[ 2 \sum_{j=1}^{f_n} \b_{j,n}^{\delta} | X_{j}^*
-Y_{j}^* |^\delta\bigg| \CG_n \Biggr]
\nonumber\\[-8pt]\\[-8pt]
&=& 2 \E\Biggl[ \sum_{j=1}^{f_n} \b_{j,n}^{\delta} \Biggr] \Delta.
\nonumber
\end{eqnarray}\looseness=0
Combining (\ref{eqwascaso1}) and (\ref{eqwascaso2}) with Proposition
\ref{Lemma2}
we obtain
\begin{eqnarray*}
&&\W_\delta^{\max(\delta,1)} \bigl(e^{-\mu(\gamma)}V_t,V_{\infty}\bigr) \\
&&\qquad\leq c
\Delta e^{-\delta\mu(\gamma)t} \sum_{n \geq0} \zeta(t,n) \E
\Biggl[ \sum
_{j=1}^{n} \b_{j,n}^{\delta}\Biggr]
\\
&&\qquad=c \Delta e^{-\delta\mu(\gamma)t} \sum_{n \geq0}
e^{-t}\bigl(1-e^{-(N-1)t}\bigr)^n \frac{({1}/({N-1}))_n}{n!}
\frac{(
{(\QQ
(\delta)+1)}/({N-1}))_n}{({1}/({N-1}))_n} \\
&&\qquad=c \Delta e^{-\delta\mu(\gamma)t} \sum_{n \geq0}
e^{-t}\bigl(1-e^{-(N-1)t}\bigr)^n \frac{(({\QQ(\delta)+1})/({N-1}))_n}{n!}
\end{eqnarray*}
with $c=1$ if $0<\gamma<\delta\leq1$ and $c=2$ if $1\leq\gamma<
\delta\leq2$,
and the thesis follows since
\[
\sum_{n \geq0} \frac{(r)_n}{n!}(1-q)^n=q^{-r}
\]
for every $q \in(0,1)$ and $r>0$.\vadjust{\goodbreak}
\end{pf*}
\begin{pf*}{Proof of Lemma~\ref{lemma-tail2}}
The proof of this lemma follows step by step the proof
of Lemma 1 in \cite{BLM2008}. %
By Lemma 9 in \cite{BLM2008}, if $\delta<\gamma/(1- \eps) $, it
suffices to prove
that
the probability distribution function of $V_\infty$, that is,
$F_\infty(x)=\int_{(-\infty,x]} \rho_{\infty,\gamma}(dy)$,
satisfies (\ref{eqFtail}) and (\ref{eqFtail2}) with the same constants
$c_0^+$ and~$c_0^-$
as the initial condition $F_0$
(possibly after diminishing $\eps$ and enlarging $K$).
The proof is based on the representation of $F_\infty$ as a mixture of
stable laws.
More precisely, let $G_\gamma$ be the distribution function
whose Fourier--Stieltjes transform is $\hat g_\gamma$ as in (\ref{chaSta});
then
\[
F_\infty(x) = \E[ G_\gamma( ( \YY)^{-1/\gamma}
x
) ],
\]
where $\YY$ has distribution $\zeta_{\gamma,\infty}$;
see Theorem~\ref{teo2}.
Since $\gamma<\delta<2\gamma$, then
there exists a finite constant $K>0$ such that
$ | 1 - c_0^+x^{-\gamma} - G_\gamma(x) | \leq K x^{-\delta}$
for $x>0$, and similarly for $x<0$;
see, for example, Sections 2.4 and 2.5 of~\cite{Zolotarev1986}.
Using that $\E[\YY]=1$ and $C:=\E[ \YY^{\delta/\gamma}]<\infty$
[by (iii) of Proposition~\ref{Pfixpoint} since $\delta
<q^*_\gamma$], it follows further that
\begin{eqnarray*}
| 1 - c_0^+ x^{-\gamma} - F_\infty(x) |
&\leq& \E[ | 1 - c_0^+ ( (\YY)^{-1/\gamma}
x)^{-\gamma} - G_\gamma(\YY^{-1/\gamma}x)| ] \\
&\leq& \E[ K (\YY)^{\delta/\gamma} x^{-\delta} ] = C K
x^{-\delta} .
\end{eqnarray*}
This proves (\ref{eqFtail}) for $F_\infty$, with $\eps=\delta
-\gamma$
and $K'=CK$.
A similar argument proves (\ref{eqFtail2}).
\end{pf*}
\begin{pf*}{Proof of Theorem~\ref{PropZolo}}
The proof follows the same line of the proof of Theorem
\ref{PropW-2}.
Assume that $\ZZ_\delta(X_0,V_\infty)<+\infty$, since otherwise the
claim is trivial.
Consider two sequences of independent and identically distributed
random variables
$(X_j)_{j \geq1}$ and $(Y_j)_{j \geq1}$, $X_j$ with common
distribution function
$F_0$ and $Y_{j}$ with the same law of $V_\infty$.
In addition assume that $(X_j)_{j \geq1},\break(Y_j)_{j \geq1}$,
$(\nu_t)_{t \geq0}$ and $(\beta_{j,n})_{j,n}$ are stochastically independent.
Recall that, as noted in the proof of Theorem
\ref{PropW-2}, $e^{\mu(\gamma)t}V_\infty$ has the same law of
$
\sum_{j=1}^{f_{\nu_t}} Y_{j} \b_{j,\nu_t}
$.
First of all it is clear, by the definition of $\ZZ_\delta$, that
\begin{eqnarray*}
\ZZ_\delta\bigl(e^{-\mu(\gamma)t} V_t,V_\infty\bigr)&=&
\ZZ_\delta\Biggl(e^{-\mu(\gamma)t}\sum_{j=1}^{f_{\nu_t}} X_{j} \b
_{j,\nu_t},
e^{-\mu(\gamma)t}\sum_{j=1}^{f_{\nu_t}}
Y_{j} \b_{j,\nu_t}\Biggr) \\
& \leq &
\sum_{n \geq0} \zeta(t,n)
\ZZ_\delta\Biggl(e^{-\mu(\gamma)t}\sum_{j=1}^{f_{n}} X_{j} \b
_{j,n},e^{-\mu(\gamma)t}\sum_{j=1}^{f_{n}}
Y_{j} \b_{j,n} \Biggr).
\end{eqnarray*}
An important property of the Zolotarev's metric $\ZZ_\delta$ is that
it is ideal of order $\delta$, that is,
\[
\ZZ_\delta(cX,cY)=c^\delta\ZZ_\delta(X,Y)
\]
(see, e.g., Theorem 1.4.2 in \cite{Zolometric}),
which yields that
\[
\ZZ_\delta\bigl(e^{-\mu(\gamma)t} V_t,
V_\infty\bigr) \leq
\sum_{n \geq0} \zeta(t,n)e^{-\delta\mu(\gamma)t}
\ZZ_\delta\Biggl(\sum_{j=1}^{f_{n}} X_{j} \b_{j,n},\sum_{j=1}^{f_{n}}
Y_{j} \b_{j,n} \Biggr).
\]
Now, by Proposition 1 in \cite{RachevRuschen},
\[
\ZZ_\delta
\Biggl(\sum_{j=1}^{f_{n}} X_{j} \b_{j,n},\sum_{j=1}^{f_{n}}
Y_{j} \b_{j,n} \Biggr) \leq\E\Biggl[\sum_{j=1}^{f_{n}} \b
_{j,n}^\delta
\Biggr]
\ZZ_\delta(X_0,V_\infty).
\]
In conclusion, we get
\[
\ZZ_\delta\bigl(e^{-\mu(\gamma)t} V_t,
V_\infty\bigr) \leq\ZZ_\delta(X_0,V_\infty) e^{-\delta\mu(\gamma
)t}\sum_{n \geq0} \zeta(t,n)
\E\Biggl[\sum_{j=1}^{f_{n}} \b_{j,n}^\delta\Biggr] .
\]
At this stage the first part of the thesis follows exactly as in the
last part of the proof of
Theorem~\ref{PropW-2}. It remains to show that, if $\gamma=2$,
$\delta
\leq3$ and $\E|X_0|^\delta<+\infty$,
then
\[
\ZZ_\delta(X_0,V_\infty) \leq\frac{1}{\Gamma(1+\delta)}[\E
|X_0|^\delta+ \E|V_\infty|^\delta],
\]
which, by Theorem~\ref{teo2}(iii) is finite.
To prove the last inequality recall that, given two random variables
$X$ and $Y$, if $2< \delta\leq3$,
then
\[
\ZZ_\delta(X,Y) \leq\frac{1}{\Gamma(1+\delta)} [\E
|X|^\delta+
\E|Y|^\delta],
\]
provided that $\E[X]=\E[Y]$ and $\E[X^2]=\E[Y^2]$; see Theorem
1.5.7 in
\cite{Zolometric}.
In our case by Theorem~\ref{teo2}(i)--(iii), $\E[X_0]=\E[V_\infty]=0$,
$\E[X_0^2]=\s_0^2$,
$\E[V^2_\infty]= \s_0^2 \int_{\RE^+} z \zeta_{\infty,2}(dz)=\s_0^2$.
\end{pf*}

\section*{Acknowledgments}

The authors would like to thank an anonymous referee for his or her
comments and careful reading of the paper. They also thank I. Gamba for
an interesting discussion.



\printaddresses

\end{document}